\documentclass[11pt,a4paper]{article}

\usepackage{graphicx,latexsym,euscript,makeidx,color,bm}
\usepackage{amsmath,amsfonts,amssymb,amsthm,thmtools,mathtools,mathrsfs,enumerate}
\usepackage[colorlinks,linkcolor=blue,citecolor=red]{hyperref}
\usepackage{tikz}

\usepackage{geometry}
\allowdisplaybreaks[4]
\def\sqr#1#2{{\vcenter{\vbox{\hrule height.#2pt
				\hbox{\vrule width.#2pt height#1pt \kern#1pt \vrule width.#2pt}
				\hrule height.#2pt}}}}
 \def\sD{\mathscr{D}}  
\def\dbE{\mathbb{E}}   
\def\dbF{\mathbb{F}} \def\sF{\mathscr{F}}  \def\cF{{\cal F}}
   
\def\dbH{\mathbb{H}}

   \def\cM{{\cal M}}

\def\dbP{\mathbb{P}}   
   
\def\dbR{\mathbb{R}}   
\def\dbS{\mathbb{S}}   
   
 \def\sU{\mathscr{U}}

\def\ss{\smallskip}             \def\hb{\hbox}
\def\ms{\medskip}              
          \def\lan{\langle}

\def\ran{\rangle}       
\def\q{\quad}               
\def\qq{\qquad}             
\def\no{\noindent}         \def\sc{\scriptscriptstyle}
         \def\scT{\sc T}
    \def\Blan{\Big\lan\!\!} \def\ds{\displaystyle}
\def\rf{\eqref}       \def\Bran{\!\!\Big\ran} 
\def\cd{\cdot}        \def\({\Big(}           
       \def\){\Big)}           \def\les{\leqslant}
\def\wt{\widetilde}   \def\[{\Big[}           \def\ges{\geqslant}
  \def\]{\Big]}           
\def\a{\alpha}       \def\l{\lambda}    \def\D{\varDelta}   
\def\b{\beta}                    
\def\d{\delta}              
\def\e{\varepsilon}      \def\L{\varLambda}
\def\f{\varphi}           
       \def\si{\sigma}    \def\Si{\varSigma}
\def\i{\infty}       \def\z{\zeta}      \def\Th{\varTheta}
\def\be{\begin{equation}}
\def\bel{\begin{equation}\label}
\def\ee{\end{equation}}
\def\bea{\begin{eqnarray}}
\def\eea{\end{eqnarray}}
\def\bt{\begin{theorem}}
\def\et{\end{theorem}}
\def\bc{\begin{corollary}}
\def\ec{\end{corollary}}
\def\bl{\begin{lemma}}
\def\el{\end{lemma}}
\def\bp{\begin{proposition}}
\def\ep{\end{proposition}}
\def\br{\begin{remark}}
\def\er{\end{remark}}
\def\ba{\begin{array}}
\def\ea{\end{array}}
\def\bd{\begin{definition}}
\def\ed{\end{definition}}
\def\5n{\negthinspace \negthinspace \negthinspace \negthinspace \negthinspace }
\def\4n{\negthinspace \negthinspace \negthinspace \negthinspace }
\def\3n{\negthinspace \negthinspace \negthinspace }
\def\2n{\negthinspace \negthinspace }
\def\1n{\negthinspace }

\newtheoremstyle{thry}
{}      
{}      
{\sl}   
{}      
{\bf}   
{.}     
{.5em}  
{}      
\theoremstyle{thry}

\newtheorem{theorem}{Theorem}[section]
\newtheorem{proposition}[theorem]{Proposition}
\newtheorem{corollary}[theorem]{Corollary}
\newtheorem{lemma}[theorem]{Lemma}

\theoremstyle{definition}
\newtheorem{definition}[theorem]{Definition}
\newtheorem{example}[theorem]{Example}

\theoremstyle{remark}
\newtheorem{remark}[theorem]{Remark}

\def\punct{}
\newtheoremstyle{dotless}{}{}{\rm}{}{\bf}{\punct}{.5em}{}
\theoremstyle{dotless}

\makeatletter
   
   \@addtoreset{equation}{section}
   \newcommand{\setword}[2]{%
   \phantomsection
   #1\def\@currentlabel{\unexpanded{#1}}\label{#2}%
   }
\makeatother

\def\3n{\negthinspace \negthinspace \negthinspace }
\def\2n{\negthinspace \negthinspace }
\def\1n{\negthinspace }
\def\bel{\begin{equation}\label}

\def\dbE{\mathbb{E}}
\def\dbF{\mathbb{F}}

\def\dbH{\mathbb{H}}

\def\dbP{\mathbb{P}}

\def\dbR{\mathbb{R}}
\def\dbS{\mathbb{S}}

\def\sD{\mathscr{D}}

\def\sF{\mathscr{F}}

\def\sU{\mathscr{U}}

\def\ds{\displaystyle}

\def\ns{\noalign{\ss}}
\def\a{\alpha}
\def\b{\beta }

\def\d{\delta}
\def\e{\varepsilon}
\def\z{\zeta}

\def\l{\lambda}

\def\si{\sigma}

\def\f{\varphi}

\def\i{\infty}
\def\D{\Delta}
\def\Th{\Theta}
\def\L{\Lambda}
\def\Si{\Sigma}

\def\O{\Omega}
\def\cF{{\cal F}}

\def\cM{{\cal M}}

\def\BS{{\bf S}}

\def\D{\Delta}
\def\Th{\Theta}
\def\L{\Lambda}
\def\Si{\Sigma}

\def\O{\Omega}

\def\no{\noindent}

\def\ss{\smallskip}
\def\ms{\medskip}

\def\q{\quad}
\def\qq{\qquad}
\def\hb{\hbox}

\def\h1{\outline{$1$}}

\def\hSi{\outline{$\Sigma$}}

\def\hTh{\outline{$\Theta$}}

\def\hh1{\outline{$1$}}
\def\hh2{\outline{$2$}}
\def\hh3{\outline{$3$}}
\def\hh4{\outline{$4$}}
\def\hh5{\outline{$5$}}
\def\hh6{\outline{$6$}}
\def\hh7{\outline{$7$}}
\def\hh8{\outline{$8$}}
\def\hh9{\outline{$9$}}
\def\hh0{\outline{$0$}}

\def\limsup{\mathop{\overline{\rm lim}}}
\def\liminf{\mathop{\underline{\rm lim}}}

\def\lan{{\langle}}
\def\ran{{\rangle}}

\def\h{\widehat}

\def\wh{\widehat}
\def\wt{\widetilde}

\def\cd{\cdot}
\def\cds{\cdots}

\def\les{\leqslant}
\def\ges{\geqslant}

\def\({\Big (}
\def\){\Big )}
\def\[{\Big[}
\def\]{\Big]}
\def\lan{\langle}
\def\ran{\rangle}

\def\bde{\begin{definition}\label}
	\def\ede{\end{definition}}
\def\bel{\begin{equation}\label}
		\def\ee{\end{equation}}
	\def\bt{\begin{theorem}\label}
		\def\et{\end{theorem}}
	\def\bc{\begin{corollary}\label}
		\def\ec{\end{corollary}}
	\def\bl{\begin{lemma}\label}
		\def\el{\end{lemma}}
	\def\bp{\begin{proposition}\label}
		\def\ep{\end{proposition}}
	\def\bex{\begin{example}\label}
		\def\ex{\end{example}}
	\def\bas{\begin{assumption}}
		\def\eas{\end{assumption}}
	\def\br{\begin{remark}\label}
		\def\er{\end{remark}}
	\def\ba{\begin{array}}
		\def\ea{\end{array}}
	\def\ed{\end{document}}

\def\rf{\eqref}

\def\square#1{\vbox{\hrule\hbox{\vrule height#1%
			\kern#1\vrule}\hrule}}
\def\rectangle#1#2{\vbox{\hrule\hbox{\vrule height#1%
			\kern#2\vrule}\hrule}}

\font\tenbb=msbm10 \font\sevenbb=msbm7 \font\fivebb=msbm5

\newfam\bbfam
\scriptscriptfont\bbfam=\fivebb \textfont\bbfam=\tenbb
\scriptfont\bbfam=\sevenbb

\makeatletter

\@addtoreset{equation}{section}
\makeatother

\usepackage{xcolor}
\makeatletter

\input pdf-trans
\newbox\qbox
\def\usecolor#1{\csname\string\color@#1\endcsname\space}
\newcommand\bordercolor[1]{\colsplit{1}{#1}}
\newcommand\fillcolor[1]{\colsplit{0}{#1}}
\newcommand\outline[1]{\leavevmode%
	\def\maltext{#1}%
	\setbox\qbox=\hbox{\maltext}%
	\boxgs{Q q 2 Tr \thickness\space w \fillcol\space \bordercol\space}{}%
	\copy\qbox%
}
\makeatother
\newcommand\colsplit[2]{\colorlet{tmpcolor}{#2}\edef\tmp{\usecolor{tmpcolor}}%
	\def\tmpB{}\expandafter\colsplithelp\tmp\relax%
	\ifnum0=#1\relax\edef\fillcol{\tmpB}\else\edef\bordercol{\tmpC}\fi}
\def\colsplithelp#1#2 #3\relax{%
	\edef\tmpB{\tmpB#1#2 }%
	\ifnum `#1>`9\relax\def\tmpC{#3}\else\colsplithelp#3\relax\fi
}
\bordercolor{black}
\fillcolor{white}
\def\thickness{.3}

\def\hSi{\outline{$\Sigma$}}

\def\hTh{\outline{$\Theta$}}

\def\h1{\outline{$1$}}

\begin{document}

\title{\bf Turnpike Property of a Linear-Quadratic Optimal Control Problem in Large Horizons with Regime Switching II: Non-Homogeneous Cases}
		
		\author{Hongwei Mei\footnote{ Department of Mathematics and Statistics, Texas Tech University, Lubbock, TX 79409, USA; email: {\tt hongwei.mei@ttu.edu}. This author is partially supported by Simons Travel Grant MP-TSM-00002835.},~~~Rui Wang\footnote{ Department of Mathematics and Statistics, Texas Tech University, Lubbock, TX 79409, USA; email: {\tt rui-math.wang@ttu.edu}},~~~\text{and}~~~
			Jiongmin Yong\footnote{Department of Mathematics, University of Central Florida, Orlando, FL 32816, USA; email: {\tt jiongmin.yong@ucf.edu}. This author is supported in part by NSF Grant DMS-2305475.}  }
				
\maketitle

\no{\bf Abstract.}
This paper is concerned with an optimal control problem for a nonhomogeneous linear stochastic differential equation having regime switching with a quadratic functional in the large time horizon. This is a continuation of the paper \cite{Mei-Wang-Yong-2025}, in which the strong turnpike property was established for homogeneous linear systems with purely quadratic cost functionals. We extend the results to the current situation. It turns out that some of the results are new even for the cases without regime switchings.

\ms
\no{\bf Keywords.}
Turnpike property, nonhomogeneous system with regime switching, stochastic optimal control, linear-quadratic problems,
stabilizability, Riccati equation.

\ms
\no{\bf AMS 2020 Mathematics Subject Classification.}  49N10, 93D20, 93E15, 93E20.

\section{Introduction}\label{Sec:Intro}

Let $(\Omega,\sF,\dbP)$ be a complete probability space on which a standard one-dimensional Brownian motion $W=\{W(t);\,t\ges 0\}$ and a  Markov chain $\a(\cd)$ with a finite state space $\cM=\{1,2,3,\cds,m_0\}$ are defined, for which they are assumed to be independent. The generator of $\a(\cd)$ is denoted by $(\l_{\imath\jmath})_{m_0\times m_0}$ (see below for details). We now denote by $\dbF^W=\{\sF^W_{t}\}_{t\ges0}$ (resp. $\dbF^\a=\{\sF_{t}^\a\}_{t\ges0}$, $\dbF=\{\sF_{t}\}_{t\ges0})$ the usual augmentation of the natural filtration generated by $W(\cd)$ (resp. by $\a(\cd)$, and by $(W(\cd),\a(\cd))$).
Consider the following {\it state equation} which is a controlled linear stochastic differential equation (SDE, for short), with regime switchings:
\bel{state}\left\{\2n\ba{ll}
\ns\ds dX(t)=\big[A(\a(t))X(t)+B(\a(t))u(t)+b(t)\big]ds\\
\ns\ds\qq\qq+\big[C(\a(t))X(t)+D(\a(t))u(t)+\si(t)\big]dW(t),\q t\in[0,T],\\
\ns\ds X(0)=x,\q\a(0)=\imath,\ea\right.\ee
%
For the coefficients of the state equation \rf{state}, we adopt the following basic assumptions:

\ms

{\bf(H1)} Let $A,C:\cM\to\dbR^{n\times n}$ and $B,D:\cM\to\dbR^{m\times n}$ be measurable.

\ms

{\bf(H2)} Let $b(\cd),\si(\cd)\in L^2_\dbF(0,T;\dbR^n)$.

\ms

Here, for any Euclidean space $\dbH$ (such as $\dbR^n,\dbR^{n\times m}$, etc.),
$$\ba{ll}
\ns\ds L^2_\dbF(0,T;\dbH)=\Big\{\f:[0,T]\times\O\to\dbH\bigm|\f(\cd)\hb{ is $\dbF$-progressively measurable, }\\
\ns\ds\qq\qq\qq\qq\qq\qq\qq\qq\qq\qq\dbE\int_0^T|\f(t)|_\dbH^2dt
<\i\Big\},\ea$$
and
$$L^{2,loc}_\dbF(0,\i;\dbH)=\bigcap_{T>0}L^2_\dbF(0,T;\dbH).$$
Also, since $\cM$ is finite, $A(\cd),B(\cd),C(\cd),D(\cd)$ are automatically bounded.

In \rf{state}, any $(x,\imath)\in\dbR^n\times\cM\equiv\sD$ is called an {\it initial pair}, $u(\cd)$, called a {\it control}, is selected from the space
$$\sU[0,T]=L_\dbF^2(0,T;\dbR^m).$$
It is well-known that for each $(x,\imath)\in\sD$ and $u(\cd)\in\sU[0,T]$, under (H1) and (H2), \rf{state} admits a unique solution $X(\cd)\equiv X(\cd\,;x,\imath;u(\cd))$, called the {\it state process}. Clearly, $X(\cd)\in L_{\dbF}^2(0,T;\dbR^n)$.
%

To measure the performance of a control $u(\cd)\in\sU[0,T]$, we introduce the following cost functional
\bel{cost[t,T]}J_{\scT}(x,\imath;u(\cd))=\dbE\(\int_0^Tg(t,
X(t),\a(t),u(t))dt\),\ee
where
\bel{g}g(t,x,\imath,u)={1\over2}\Blan\begin{pmatrix}Q(\imath)& S(\imath)^\top \\ S(\imath)&R(\imath)\end{pmatrix}\begin{pmatrix}x\\ u\end{pmatrix},\begin{pmatrix}x\\ u\end{pmatrix}\Bran
+\Blan\begin{pmatrix}q(t) \\ r(t)\end{pmatrix},\begin{pmatrix}x\\ u\end{pmatrix}\Bran,\ee
with $Q(\cd),S(\cd),R(\cd)$ being suitable matrix valued maps and some stochastic processes $q(\cd),r(\cd)$. Here, the superscript $\top$ denotes the transpose of matrices; $\lan\cd\,,\cd\ran$ denotes the inner product of two vectors (possibly in different spaces). In what follows, we denote $\dbS^n$, $\dbS^n_+$ and $\dbS^n_{++}$ to be the sets of all $(n\times n)$ symmetric, positive semi-definite, and positive definite matrices, respectively. For the weights in the cost functional \rf{cost[t,T]}, we adopted the following assumption.

\ms

{\bf(H3)} Suppose that $Q(\imath)\in\dbS^n_{++}$, $R(\imath)\in\dbS^m_{++}$, $S(\imath)\in\dbR^{n\times m}$ such that
\bel{Q}Q(\imath)-S(\imath)^\top R(\imath)^{-1}S(\imath)\in\dbS^n_{++}.\ee

{\bf(H4)} Let $q(\cd)\in L^2_\dbF(0,T;\dbR^n)$ and $r(\cd)\in L^2_\dbF(0,T;\dbR^m)$.

\ms

Now it is natural to consider the following optimal control problem, under (H1)--(H4).

{\bf Problem (LQ)$_{\scT}$.} For a given initial pair $(x,\imath)\in\sD$, find a control $\bar u_{\scT}^{x,\imath}(\cd)\in\sU[0,T]$ such that
\bel{opt-u3}J_{\scT}(x,\imath;\bar u_{\scT}^{x,\imath}(\cd))=\inf_{u(\cd)\in\sU[0,T]}J_{\scT}(x,\imath;u(\cd))\equiv V_{\scT}(x,\imath).\ee

The above problem is referred to as a (nonhomogeneous) {\it linear-quadratic} (LQ, for short) {\it optimal control problem} over a finite horizon with regime switchings (see \cite{Zhang-Li-Xiong-2021, Mei-Wei-Yong-2025} for examples). We call $\bar u^{x,\imath}_{\scT}(\cd)$ an {\it open-loop optimal control process}, $\bar X^{x,\imath}_{\scT}(\cd)$ the corresponding {\it open-loop optimal state process}, and $(\bar X^{x,\imath}_{\scT}(\cd),\bar u^{x,\imath}_{\scT}(\cd))$  the {\it open-loop optimal pair} of Problem (LQ)$_{\scT}$, respectively. In addition, we call $V_{\scT}(x,\imath)$ the {\it value function} of Problem (LQ)$_{\scT}$.


Under some general mild assumptions (which will be given in later sections), it can be proved that Problem (LQ)$_{\scT}$ admits a unique open-loop optimal control $\bar u_{\scT}^{x,\imath}(\cd)\in\sU[0,T]$ with the optimal state process $\bar X_{\scT}^{x,\imath}(\cd)\equiv X(\cd\,;x,\imath;\bar u_{\scT}^{x,\imath}(\cd))\in L^2_\dbF(0,T;\dbR^n)$. For the cases without regime switchings, people found that, under the so-called stabilizability condition (see below), for some stochastic processes $(\bar X_{\sc\i}(\cd),\bar u_{\sc\i}(\cd))$, and some constants $\b,K>0$, all are independent of $0<T<\i$, (in what follows, $K>0$ will be a generic constant which can be different from line to line) such that
\bel{STP0}\dbE\big(|\bar X_{\scT}^{x,\imath}(t)\1n-\1n\bar X_{\sc\i}(t)|^2\2n+\1n|\bar u_{\scT}^{x,\imath}(t)\1n-\1n\bar u_{\sc\i}(t)|^2\big)\les
\1n K(e^{-\b t}\1n+\1n e^{-\b(T-t)}),\q\forall t\in[0,T].\ee
Such an asymptotic behavior of the optimal pair $(\bar X_{\scT}^{x,\imath}(\cd),\bar u_{\scT}^{x,\imath}(\cd))$ as $T\to\i$ is called the {\it strong turnpike property} (STP, for short) of Problem (LQ)$_{\scT}$. The main feature of \rf{STP0} is that the (open-loop) optimal pair $(\bar X^{x,\imath}_{\scT}(\cd),\bar u_{\scT}^{x,\imath}(\cd))$ will be very close to a $T$-independent pair $(\bar X_{\sc\i}(\cd),\bar u_{\sc\i}(\cd))$ for all $t$ in the middle range of $[0,T]$ (i.e., $t\in[\e T,(1-\e)T]$ for some $\e\in(0,{1\over2})$.)

Research on turnpike phenomenon was begun by Ramsey (\cite{Ramsey-1928}) in 1928, followed by von Neumann (\cite{Neumann-1945}) in 1945, and  Dorfman--Samuelson--Solow (\cite{Dorfman-Samuelson-Solow-1958}) in 1958, who coined the name.  Since then, the turnpike property has been found to hold for a large class of (deterministic, finite or infinite dimensional) optimal control problems. Numerous relevant results can be found in \cite{McKenzie-1976,Carlson-Haurie-Leizarowitz-1991,Damm-Grune-Stieler-Worthmann-2014,
Trelat-Zuazua-2015,Zuazua-2017,Grune-Guglielmi-2018,Zaslavski-2019,
Lou-Wang-2019,Breiten-Pfeiffer-2020,Sakamoto-Zuazua-2021,Faulwasser-Grune-2022} and the references cited therein. Since beginning of 1970, several authors studied the problem from portfolio aspect showing that for certain maximization problems of the utility for investments, the turnpike properties were established, mainly under proper assumptions on the utility functions (see \cite{Leland-1972, Ross-1974,Humberman-Ross-1983,Cox-Huang-1992,Jin-1998,
Huang-Zariphopoulou-1999,Dybvig-Rogers-Back-1999,Bian-Zheng-2019,
Geng-Zariphopoulou-2025}). Recently, a systematic investigation for continuous-time stochastic optimal LQ control problems was begun by the work of Sun--Wang--Yong in the early of 2020 (\cite{Sun-Wang-Yong-2022}), followed by the works  \cite{Conforti-2023,Chen-Luo-2023,Sun-Yong-2024S,Sun-Yong-2024J,Jian-Jin-Song-Yong-2024,Schiessl-Baumann-Faulwasser-Grune-2024,
Bayraktar-Jian-2025}. In particular, turnpike property for stochastic LQ control problems with regime switching has been studied by the authors in \cite{Mei-Wang-Yong-2025} when the linear SDE is homogeneous and the cost functional is purely quadratic. Naturally, one may ask if the results of \cite{Mei-Wang-Yong-2025} are true for nonhomogeneous problems, with the cost functional also having linear terms. The purpose of the current paper is to give a positive answer to this question, with additional techniques.

More precisely, combing those in \cite{Mei-Wang-Yong-2025}, with some additional assumptions (see below), we will refine \rf{STP0} as follows: there exists a function $h(\cd)\ges0$ and constants $\b,K>0$, all are independent of $0<T<\i$, such that the following refined STP holds:
\bel{refinedSTP}\ba{ll}
\ns\ds\dbE\[|\bar X_{\scT}^{x,\imath}(t)-\bar X^{x_{\sc\i},\imath}_{\sc\i}(t)|^2+\int_0^t|\bar u_{\scT}^{x,\imath}(s)-\bar u^{x_{\sc\i},\imath}_{\sc\i}(s)|^2ds\]\\
\ns\ds\q\les K\[e^{-\b t}|  x-x_{\sc\i}|^2+ e^{-\b(T-t)}\(e^{-\b t}|x|^2+h(t)\)\],\qq\forall t\in[0,T],\ea\ee
with $(x,\imath),(x_\i,\imath)\in\sD$ being two possibly different initial pairs.

In particular, if we strengthen (H2) and (H4) to the following:

\ms

{\bf(H2)$'$} Let $b(\cd),\si(\cd),q(\cd)\in L^2_{\dbF^\a}(0,T;\dbR^n)$ and $r(\cd)\in L^2_{\dbF^\a}(0,T;\dbR^m)$ be bounded for any $T>0$.

\ms

\no Then the above \rf{refinedSTP} can be strengthened to the following:
\bel{refinedSTP*}\ba{ll}
\ns\ds\dbE\[|\bar X_{\scT}^{x,\imath}(t)-\bar X^{x_{\sc\i},\imath}_{\sc\i}(t)|^2+|\bar u_{\scT}^{x,\imath}(t)-\bar u^{x_{\sc\i},\imath}_{\sc\i}(t)|^2\]\\
\ns\ds\q\les K\[e^{-\b t}|  x-x_{\sc\i}|^2+ e^{-\b(T-t)}\(e^{-\b t}|x|^2+h(t)\)\],\qq\forall t\in[0,T],\ea\ee

Now, we indicate three types of asymptotic behaviors of the open-loop optimal pair to the relevant Problem (LQ)$_{\scT}$:

$\bullet$ {\bf Homogeneous Case}: Let $b(\cd),\si(\cd),q(\cd),r(\cd)$ be all 0, and (H1), (H3) hold. In this case, $h(t)\equiv0$. This case, fully treated in \cite{Mei-Wang-Yong-2025}, is singled out since it catches one of the most essential features of the problem: the convergence of the solutions to differential Riccati equations (DREs, for short) to that of the algebraic Riccati equation (ARE, for short).

$\bullet$ {\bf Integrable Case}: $b(\cd),\si(\cd),q(\cd)\in L_\dbF^2(0,\i;\dbR^n)$ and $r(\cd)\in L_\dbF^2(0,\i;\dbR^m)$. In this case, $h(\cd)$ is a non-negative integrable function on $[0,\i)$. Due to the appearance of the nonhomogeneous (square integrable) terms $b(\cd),\si(\cd)$ in the state equation \rf{state} and linear (square integrable) weights $q(\cd),r(\cd)$ in the cost functional \rf{cost[t,T]}, some backward stochastic differential equations (BSDEs, for short) will be involved. From this case, we can see how far one can go (by this approach). This (even for the cases without switching) is new in the literature, since those non-homogeneous terms were assumed to be constants (\cite{Sun-Wang-Yong-2022, Sun-Yong-2024S}) or periodic (\cite{Sun-Yong-2024J}).

$\bullet$ {\bf Local-Integrable Case}: For any $0<T<\i$, $b(\cd),\si(\cd),q(\cd)\in L_\dbF^2(0,T;\dbR^n)$ and $r(\cd)\in L_\dbF^2(0,T;\dbR^m)$ with some additional assumptions. In this case, we can take $h(t)\equiv1$. This is new again, even for the optimal control problems without regime switchings. Without doubt, the LQ ergodic control problem will be involved. We can see that the previous results in \cite{Sun-Wang-Yong-2022,Sun-Yong-2024S,Sun-Yong-2024J} are some special cases of that in the current paper.

The rest of the paper is arranged as follows. In Section \ref{sec:ocp}, we recall results for Problem (LQ)$_{\scT}$, with $0<T<\i$. Section \ref{sec:sta} is to devote to some asymptotic behavior of the open-loop optimal pair $(\bar X_{\scT}^{x,\imath}(\cd),\bar u_{\scT}^{x,\imath}(\cd))$ as $T\rightarrow\i$, under the stabilizability condition, including the identification of the limit pair $(\bar X_{\sc\i}^{x,\imath}(\cd),\bar u_{\sc\i}^{x,\imath}(\cd))$. 
Then our main results on STP are proved in Section \ref{sec:STP}. In particular, we verify the optimality of $(\bar X_{\sc\i}^{x,\imath}(\cd),\bar u_{\sc\i}^{x,\imath}(\cd))$ in two different cases: integrable and local-integrable cases in Section \ref{sec:opt}.  Then some concluding remarks are made in Section \ref{sec:con}. Finally, some proofs are relegated in Section \ref{sec:pro}.

\section{Optimal Control of Problem (LQ)$_{\scT}$} \label{sec:ocp}

In this section, we will recall the optimal control and its closed-loop representation for Problem (LQ)$_{\scT}$ with $0<T<\i$. Our results are some special cases of those found in \cite{Zhang-Li-Xiong-2021}.
\ms

Recall that $\a(\cd)$ is a Markov chain whose state space $\cM$ is finite. Thus, we may let its generator be $(\l_{\imath\jmath})_{m_0\times m_0}\in\dbR^{m_0\times m_0}$, which is a real matrix so that the following hold:
\bel{q-prop}\l_{\imath\jmath}>0,\q\imath\ne\jmath;\qq\sum_{\jmath=1}^{m_0}
\l_{\imath\jmath}=0,\q\imath\in\cM.\ee
We now proceed with a {\it martingale measure of Markov chain} $\a(\cd)$. For $\imath\ne\jmath$, we define
$$\ba{ll}
\ns\ds\wt M_{\imath\jmath}(t):=\sum_{0\leq s\les t}{\bf1}_{[\a(s_-)=\imath]}{\bf1}_{[\a(s)=\jmath]}\equiv\hb{accumulative jump number from $\imath$ to $\jmath$ 
in $(0,t]$},\\
\ns\ds\lan\wt M_{\imath\jmath}\ran(t):=\int_0^t
\l_{\imath\jmath}{\bf1}_{[\a(s-)=\imath]}ds,
\q M_{\imath\jmath}(t):=\wt M_{\imath\jmath}(t)-\lan\wt M_{\imath\jmath}\ran(t),\qq s\ges0.\ea$$
The above $M_{\imath\jmath}(\cd)$ is a square-integrable martingale (with respect to $\dbF^\a$). For convenience, we let
$$M_{\imath\imath}(t)=\wt M_{\imath\imath}(t)=\lan\wt M_{\imath\imath}\ran(t)=0,\qq s\ges0.$$
Then $\{M_{\imath\jmath}(\cd)\bigm|\imath,\jmath\in\cM\}$ is the {\it martingale measure} of Markov chain $\a(\cd)$. If $\dbH$ is a Euclidean space and $F:\cM\to\dbH$ is measurable, then
\bel{dF}d[F(\a(t))]=\L[F](\a(t))ds+\sum_{\imath,\jmath\in\cM}[F(\jmath)
-F(\imath)]{\bf1}_{\{\a(s^-)=\imath\}}dM_{\imath\jmath}(t),\ee
where (see \rf{q-prop})
\bel{L[F]}\L[F](\imath)=\sum_{\jmath\ne\imath}\l_{\imath\jmath}
F(\jmath)\equiv\sum_{\imath,\jmath\in\cM}\l_{\imath\jmath}[F(\jmath)-F(\imath)].\ee
This is a special case of \cite{Yin-Zhang-2010}, Section 2.7, or \cite{Yin-Zhu-2010}, Section 2.2. In fact,
\begin{align*}&F(\a(t))-F(\a(0))\\
&=\sum_{0\les s\les t}[F(\a(s))-F(\a(s^-))]=\sum_{0\les r\les s}\sum_{\imath,\jmath\in\cM}[F(\jmath)-F(\imath)]{\bf1}_{
\{\a(s)=\jmath\},\a(s^-)=\imath\}}\\
&=\int_0^t\sum_{\imath,\jmath\in\cM}[F(\jmath)-F(\imath)]
{\bf1}_{\{\a(s^-)=\imath\}}d\wt M_{\imath\jmath}(s)=\int_0^t\sum_{\imath,\jmath\in\cM}\l_{\imath\jmath}[F(\jmath)-F(\imath)]
{\bf1}_{\{\a(s^-)=\imath\}}ds\\
&+\int_0^t\sum_{\imath,\jmath\in\cM}[F(\jmath)-F(\imath)]
{\bf1}_{\{\a(s^-)=\imath\}}d\(\wt M_{\imath\jmath}(s)-\l_{\imath\jmath}s\)\\
&=\int_0^t\L[F](\a(s))ds+\int_0^s\sum_{\imath,\jmath\in\cM}[F(\jmath)-F(\imath)]
{\bf1}_{\{\a(s^-)=\imath\}}dM_{\imath\jmath}(s).\end{align*}
Thus, we have \rf{dF}.

Now, let $\dbF_-$ be the smallest filtration containing $\{\cF_t^W\}_{t\ges 0}$ and $\{\cF_{t-}^\a\}_{t\ges 0}$ augumented with all $\dbP$-null sets.
To define the stochastic integral with respect to such a martingale measure, we need to introduce the following Hilbert spaces
$$\ba{ll}
\ns\ds M^2_{\dbF_-}(t,T;\dbH)=\Big\{\f(\cd\,,\cd)=(\f(\cd\,,1),\cds,\f(\cd\,,m_0))\bigm|
\f(\cd\,,\cd)\hb{ is $\dbH$-valued and $\dbF_-$-measurable }\\
\ns\ds\qq\qq\qq\qq\qq \hb{ with } \dbE\int_t^T\sum_{\imath\neq\jmath}|
\f(s,\jmath)|^2\l_{\imath\jmath}{\bf1}_{[\a(s)=\imath]}ds<
\i,\q\forall\imath,\jmath\in\cM
\Big\}.\ea$$
Now, for any $\f(\cd)\in M^2_{\dbF_-}(t,T;\dbH)$, we define its stochastic integral against $dM$ by the following:
$$\int_t^T\f(s)dM(s):=\sum_{\jmath\ne\imath}\int_{[t,T]}\f(r,\jmath)
{\bf1}_{[\a(s^-)=\imath]}dM_{\imath\jmath}(s),$$
whose quadratic variation is
$$\dbE\(\int_t^T\f(s)dM(s)\)^2=\dbE\int_t^T\sum_{\imath\neq\jmath}|
\f(s,\jmath)|^2\l_{\imath\jmath}{\bf1}_{[\a(s)=\imath]}ds.
$$

Now we state the following results concerning Problem (LQ)$_{\scT}$, whose proof can be found in \cite{Zhang-Li-Xiong-2021,Mei-Wei-Yong-2024} (see also \cite{Sun-Yong-2020}).

\bp{propBSDET} \sl Let {\rm (H1)--(H4)} hold.

{\rm(i)} For each $\imath\in\cM$, the following DRE  admits a unique  uniformly regular solution $P_{\scT}(\cd,\imath)\in C(0,T;\dbS_{++}^n)$ ($\imath\in\cM$):
\bel{dsE}\left\{\2n\ba{ll}
\ds\dot P_{\scT}+\L[P_{\scT}]+P_{\scT}A+A^\top P_{\scT}+C^\top P_{\scT}C+Q\\
\ns\ds-(P_{\scT}B\1n+\1n C^\top\1n P_{\scT}D\1n+\1n S^\top)(R\1n+\1n D^\top\1n P_{\scT}D)^{-1}(B^\top\1n P_{\scT}\1n+\1n D^\top\1n P_{\scT}C\1n+\1n S)\1n=\1n0,\q t\in[0,T],\\
\ns\ds P_{\scT}(T)=0,\ea\right.\ee
i.e., it is a solution of \rf{dsE} and for some $T$-independent constant $\d>0$, it holds
\bel{R+DPD>0}\wt R_{\scT}(t,\imath)\equiv R(\imath)+D^\top(\imath)P_{\scT}(t,\imath)D(\imath)\ges\d I,\q\forall(t,\imath)\in[0,T]\times\cM.\ee

{\rm(ii)} There exists a unique adapted solution  $(\eta_{\scT}(\cd),\z_{\scT}(\cd),\z^M_{\scT}(\cd))
\in L_{\dbF}^2(0,T;\dbR^n)\times L_{\dbF}^2(0,T;\dbR^n)\times M_{\dbF_-}^2(0,T;\dbR^n)$ solving the following BSDE on $[0,T]$:
\bel{BSDET}\left\{\2n\ba{ll}
\ds d\eta_{\scT}(t)=-\Big(A^{\Th_T}(t,\a(t))^\top \eta_{\scT}(t)+C^{\Th_T}(t,\a(t))^\top\z_{\scT}(t)+\f_{\scT}(t,\a(t))\Big)ds\\
\ns\ds\qq\qq\q+\z_{\scT}(t)dW(t)+\z^M_{\scT}(t)dM(t),\\
\ns\ds\eta_{\scT}(T)=0,\ea\right.\ee
where
\bel{A^Th_T}\ba{ll}
\ns\ds A^{\Th_T}(t,\imath)=A(\imath)+B(\imath)\Th_{\scT}(t,\imath),\qq
C^{\Th_T}(t,\imath)=C(\imath)+D(\imath)\Th_{\scT}(t,\imath),\\
\ns\ds\Th_{\scT}(t,\imath)\1n=\1n-\wt R_{\scT}(t,\imath)^{-1}[B(\imath)^\top\1n P_{\scT}(t,\imath)\1n+\1n D(\imath)^\top\1n P_{\scT}(t,\imath)C(\imath)\1n+\1n S(\imath)],\\
\ns\ds\f_{\scT}(t,\imath)=P_{\scT}(t,\imath)b(t)+C^{\Th_T}(t,\imath)^\top\1n P_{\scT}(t,\imath)\si(t)\1n+\1n\Th_{\scT}(t,\imath)^\top\1n r(t)+q(t).\ea\ee

{\rm(iii)} For each $(x,\imath)\in\sD$,  the unique open-loop optimal control $\bar u^{x,\imath}_{\scT}(\cd)\in \sU[0,T]$ admits a closed-loop representation
\bel{optT2}\bar u_{\scT}^{x,\imath}(t)=\Th_{\scT}(t,\a(t))\bar X^{x,\imath}(t)+v_{\scT}(t,\a(t)),\qq t\in[0,T],\ee
where $\bar X^{x,\imath}(\cd)$ is the corresponding optimal state process and
\bel{v_T}\ba{ll}
\ns\ds v_{\scT}(t,\imath)=-\wt R_{\scT}(t,\imath)^{-1}[D(\imath)^\top P_{\scT}(t,\imath)\si(t)+B(\imath)^\top\eta_{\scT}(t)+D(\imath)^\top
\z_{\scT}(t)+r(t)],\\
\ns\ds\qq\qq\qq\qq\qq\qq\qq\qq\qq\qq\qq(t,\imath)\in[0,T]\times\cM.\ea\ee
In this case, the optimal closed-loop system reads
\bel{closed-loop}\left\{\2n\ba{ll}
\ds d\bar X_{\scT}^{x,\imath}(t)=[A^{\Th_{\scT}}(t,\a(t))\bar X_{\scT}^{x,\imath}(t)+B(\a(t))v_{\scT}(t,\a(t))+b(t)]ds\\
\ns\ds\qq\qq+[C^{\Th_{\scT}}(t,\a(t))\bar X_{\scT}^{x,\imath}(t)+D(\a(t))v_{\scT}(t,\a(t))+\si(t)]dW(t),\\
\ns\ds\bar X_{\scT}^{x,\imath}(0)=x,\qq\a(0)=\imath.\ea\right.\ee

\ep

Having the open-loop optimal pair $(\bar X_{\scT}^{x,\imath}(\cd),\bar u_{\scT}^{x,\imath}(\cd))$ of Problem (LQ)$_{\scT}$, our next goal is to obtain the asymptotic behavior of this optimal pair as $T\to\i$. This will be carefully investigated in the following section.

\section{Asymptotic Behavior of Optimal Controls} \label{sec:sta}

In this section, we will investigate the asymptotic behavior of the open-loop optimal pair $(\bar X^{x,\imath}(\cd),\bar u^{x,\imath}(\cd))$ as $T\to\i$. A so-called stabilizability condition is required. First, we will consider the asymptotic behavior of $\Th_{\scT}(\cd)$. This part has been fully studied in \cite{Mei-Wang-Yong-2025}. For readers' convenience, we recall the main results here.
Based on this, we will further derive the asymptotic behavior for $v_{\scT}(\cd)$.

Let
$$\ba{ll}
\ns\ds\hTh=\Big\{\Th:\cM\to\dbR^{m\times n}\bigm|\Th(\cd)\hb{ is measurable}\Big\},\\
\ns\ds\hSi=\Big\{\Si:\cM\to\dbS^n_{++}\bigm|\Si(\cd)\hb{ is measurable}\Big\},\ea$$
and let us consider the following linear SDE with a regime switching governed by a Markov chain:
\bel{state3}\left\{\2n\ba{ll}
\ds dX(t)=A(\a(t))X(t)dt+C(\a(t))X(t)dW(t),\qq t\in[0,\i),\\
\ns\ds X(0)=x,\q\a(0)=\imath.\ea\right.\ee
The above system is denoted by $[A,C]$. Under (H1), such a system is well-posed. If $X(\cd)\equiv X(\cd\,;x,\imath)$ is the solution of the above corresponding to $(x,\imath)\in\sD$.
We now introduce the following definition.

\bde{stable} \rm (i) System $[A,C]$ is said to be {\it stable} if for any $(x,\imath)\in\sD$, $X(\cd\,;x,\imath)\in L_{\dbF}^2(0,\i;\dbR^n)$,

(ii) System $[A,C]$ is said to be {\it dissipative} if one could find a $\Si(\cd)\in\hSi$ and a $\d>0$ so that
\bel{dissipative}\(\L[\Si]+\Si A+A^\top\Si+C^\top\Si C\)(\imath)\les-\d\Si(\imath),\qq\imath\in\cM.\ee

\ede

The following definition is adopted from \cite{Mei-Wei-Yong-2025}.

\bde{Def-stab-1} (i) System $[A,C;B,D]$ is said to be {\it stabilizable} if one can find a map $\Th(\cd)\in\hTh$, so that for $[A^\Th,C^\Th]$ is stable, where (see \rf{A^Th_T})
\bel{A^Th}A^\Th(\imath)=A(\imath)+B(\imath)\Th(\imath),\q C^\Th(\imath)=C(\imath)+D(\imath)\Th(\imath).\ee
In this case, the map $\Th(\cd)$ is called a {\it stabilizer} of $[A,C;B,D]$. The set of all possible stabilizers of system $[A,C;B,D]$ is denoted by $\BS[A,C; B,D]$.

(ii) The map $\Th(\cd)\in\hTh$ is called a {\it dissipating strategy} of system $[A,C;B,D]$ if there exists a $\d>0$ and a $\Si(\cd)\in\hSi$ such that \rf{dissipative} holds with $[A,C]$ replaced by $[A^\Th,C^\Th]$, i.e.,
\bel{dissipative*}\(\L[\Si]+\Si A^\Th+(A^{\Th})^\top\Si+(C^\Th)^\top\Si C^\Th\)(\imath)\les-\d\Si(\imath),\qq\imath\in\cM.\ee

\ede

Since the state space $\cM$ of the Markov chain $\a(\cd)$ is finite, the following is true (see \cite{Mei-Wei-Yong-2025}, Proposition 3.7).

\bp{3.3} System $[A,C;B,D]$ is stabilizable if and only if it admits a dissipating strategy.

\ep

It is known that stabilizability of $[A,C;B,D]$ is necessary for studying LQ problems in an infinite time horizon even for the problems without regime switchings (\cite{Sun-Yong-2020, Mei-Wei-Yong-2025}). Thus, we accept the following assumption.

\ms

{\bf(H5)} System $[A,C;B,D]$ is stabilizable, i.e., $\BS[A,C;B,D]\ne\varnothing$.
\ms

To find the asymptotic behavior of $P_{\scT}(t,\imath)$ as well as $\Th_{\scT}(t,\imath)$ (as $T\rightarrow\i$), we introduce the following ARE:
\bel{ARE}\ba{ll}
\ns\ds\L[P_{\sc\i}]+P_{\sc\i}A+A^\top P{\sc\i}+C^\top P_{\sc\i}C+Q\\
\ns\ds\q-(B^\top P_{\sc\i}+D^\top P_{\sc\i}C+S)^\top(R+D^\top P_{\sc\i}D)^{-1}(B^\top P_{\sc\i}+D^\top P_{\sc\i}C+S)=0.\ea\ee
The following is the key result obtained in \cite{Mei-Wang-Yong-2025}. The main feature is the convergence.

\bp{3.4} \sl Let {\rm(H1), (H3)} and {\rm(H5)} hold. Then

\ms

{\rm (i)} ARE \rf{ARE} admits a unique regular solution $P_{\sc\i}(\cd):\cM\to\dbS^n_{++}$, i.e., it is a solution of \rf{ARE} such that
\bel{R+DPD>0*}\wt R_{\sc\i}(\imath)\equiv R(\imath)+D(\imath)^\top P_{\sc\i}(\imath)D(\imath)\ges\d I,\q\imath\in\cM,\ee
for some $\d>0$, and
\bel{3.8}\Th_{\sc\i}(\cd)=-\wt R_{\sc\i}(\cd)^{-1}[B(\cd)^\top P_{\sc\i}(\cd)+D(\cd)^\top P_{\sc\i}(\cd)C(\cd)+S(\cd)]\in\BS[A,C;B,D],\ee
i.e., there exists a $\d>0$ and $\Si_{\sc\i}(\cd)\in\hSi$ such that (by Proposition \ref{3.3})
\bel{Ly1}\(\L[\Si_{\sc\i}]+\Si_{\sc\i} A^{\Th_\i}+(A^{\Th_\i})^\top\Si_{\sc\i}+(C^{\Th_\i})^\top\Si_{\sc\i}
C^{\Th_\i}\)(\imath)\les-\d\Si_{\sc\i}(\imath),\q\forall\imath\in\cM.\ee

{\rm (ii)} For any given $t\in[0,\i)$, the following convergence holds
\bel{P to P}P_{\scT}(t,\imath)=P_{\scT-t}(0,\imath)\nearrow P_{\sc\i}(\imath),\qq\hb{as } T\nearrow\i,\q\forall\imath\in\cM.\ee
Moreover, there exists a $\d>0$ so that (for some absolute constants $K,\d>0$)
\bel{P-P}
0\les P_{\sc\i}(\imath)-P_{\scT}(t,\imath)\les Ke^{-\d(T-t)}I,\qq t\in[0,T],\ee
and consequently,
\bel{Th-Th}|\Th_{\sc\i}(\imath)-\Th_{\scT}(t,\imath)|\les Ke^{-\d(T-t)},\qq t\in[0,T],\ee

{\rm (iii)} There exists a constant $0<T_0<T$ with $T-T_0\ges 0$ large enough such that
\bel{Ly2}\ba{ll}
\ns\ds\L[\Si_{\sc\i}(\cd)](\imath)+\Si_{\sc\i}(\imath) A^{\Th_{\scT}}(t,\imath)+A^{\Th_{\scT}}(t,\imath)^\top
\Si_{\sc\i}(\imath)\\
\ns\ds\qq\qq+C^{\Th_{\scT}}(t,\imath)^\top\Si_{\sc\i}(\imath) C^{\Th_{\scT}}(t,\imath)\les-{\d\over2}\Si_{\sc\i}(\imath),\qq
t\in[0,T-T_0],\ea\ee
where $A^{\Th_T}$ and $C^{\Th_T}$ are given by \rf{A^Th_T}.
%
%

\ep

\begin{proof}. \rm (i) and (ii) are derived in \cite{Mei-Wang-Yong-2025}. (iii) is concluded from \rf{Th-Th} and inequality \rf{Ly1} since
$$\ba{ll}
\ns\ds|A^{\Th_T}(t,\imath)-A^{\Th_\i}(\imath)|+|C^{\Th_T}(t,\imath)-
C^{\Th_\i}(\imath)|\\
\ns\ds\les\big(|B(\imath)|+|D(\imath)|\big)|\Th_{\scT}(t,\imath)
-\Th_{\sc\i}(\imath)|\les Ke^{-\d(T-t)}.\ea$$
The choice of $T_0>0$ is such that
$Ke^{-\d(T-t)}\les Ke^{-\d T_0}$
is small enough for all $t\in[0,T-T_0]$. 

\end{proof}

We note that when $t$ is close to $T$, say, $0<T-t\les T_0$, the right-hand sides of \rf{P-P} and \rf{Th-Th} might not be small. In other words, only if $t$ is far away from $T$, say, $T-t\ges
T_0$, the right-hand sides of \rf{P-P} and \rf{Th-Th} will be small, and \rf{Ly2} will be true.

Now, we introduce the following assumption.

\ms

{\bf(H6)} Let $b(\cd),\si(\cd),q(\cd)\in L^{2,loc}_\dbF(0,\i;\dbR^n)$, and $r(\cd)\in L^{2,loc}_\dbF(0,\i;\dbR^m)$.

\ms



 Under the dissipativity assumption, we are able to further derive the following proposition concerning with the existence and uniqueness of the adapted solutions to BSDEs \rf{BSDET}, together with several useful
estimates. The proof of the proposition is quite lengthy and will be given in Section \ref{sec:pro}. Write
\begin{align}\label{defxi}\xi(t)\equiv\dbE[|b(t)|^2+|\si(t)|^2+|q(t)|^2+|r(t)|^2].
\end{align}

\bp{3.5} \sl Let {\rm(H1)--(H6)} hold. Let $\d>0$ given by \rf{Ly1} and $T_0$ be that in {\rm(iii)} of Proposition \ref{3.4}. Then it follows that
\begin{align}&\dbE|\eta_{\scT}(t)|^2+\dbE\int_t^T e^{-{\d\over4}(s-t)}\sum_{\jmath\neq\imath}\lambda_{\imath\jmath}|\z^M_{\scT}(s,\jmath)|^2{\bf 1}_{[\a(s)=\imath]}ds\nonumber\\
&\q+\dbE\int_t^T e^{-{\d\over4}(s-t)}|\z_{\scT}(s)|^2ds\les K\int_t^T e^{-{\d\over4}(s-t)}\xi(s)ds.\label{stabsde}
\end{align}
For any $T'>T>T_0$, it also holds that
\begin{align}
&\dbE|\eta_{\scT}(t)-\eta_{\scT'}(t)|^2+\dbE\int_t^{T}e^{-{\d\over4}(s-t)}\sum_{\jmath\neq\imath}\lambda_{\imath\jmath}|\z^M_{\scT}(s,\jmath)-\z^M_{\scT'}(s,\jmath)|^2{\bf 1}_{[\a(s)=\imath]}ds\nonumber\\
&\q+\dbE\int_t^{T}
e^{-{\d\over4}(s-t)}|\z_{\scT}(s)-\z_{\scT'}(s)|^2ds\les Ke^{-{\d\over8}(T-s)}\int_t^{\scT'} e^{-{\d\over4}(s-t)}\xi(s)ds.\label{E[eta-eta]}
\end{align}

\ep

 The above proposition, \eqref{E[eta-eta]} particularly, suggests  the following assumption.\ms

{\bf(H6)$'$} For $\d>0$ given by (iii) of Proposition \ref{3.5}, the following holds:
\bel{int<i}\sup_{t\in[0,\i)}\int_0^\i e^{-{\d\over4}|t-s|}\xi(s)ds<\i.\ee

\ms

Even though (H6)$'$ is a little stronger than (H6), we see that  (H6)$'$ holds if $\xi(\cd)$ is measurable and bounded. Thus, (H6)$'$ covers most interesting cases.

\ms

Under (H6)$'$, taking $T\to\i$ in \rf{BSDET}, formally, we have the following BSDE on $[0,\i)$:
\begin{align}\label{BSDEIN}& d\eta_{\sc\i}(t)=-\Big(A^{\Th_{\sc\i}}(t,\a(t))^\top \eta_{\sc\i}(t)+C^{\Th_{\sc\i}}(t,\a(t))^\top\z_{\sc\i}(t)+\f_{\sc\i}(t,\a(t))\Big)dt\nonumber\\
&\qq\qq\qq+\z_{\sc\i}(t)dW(t)+\z^M_{\sc\i}(t)dM(t),
\end{align}
with
\bel{f_i}\f_{\sc\i}(t,\imath)=P_{\sc\i}(t,\imath)b(t)+C^{\Th_\i}
(t,\imath)^\top\1n P_{\sc\i}(t,\imath)\si(t)\1n+\1n\Th_{\sc\i}(t,\imath)^\top\1n r(t)+q(t).\ee
A similar argument to the proof of Proposition \ref{3.5}  can show that \eqref{BSDEIN} admits a unique solution  BSDE  $$(\eta_{\i }(\cd), \z_{\i }(\cd),\z_{\i }^M(\cd))\in L_{\dbF}^{2}(0,T;\dbR^n)\times L_{\dbF}^{2}(0,T;\dbR^n)\times M_{\dbF_-}^{2,T}(0,T;\dbR^n),$$
for any $T>0$. Moreover, \eqref{stabsde} holds for $T=\infty$  and  \eqref{E[eta-eta]} holds for any $T>T_0$ and  $T'=\i.$
\begin{remark}\rm
(1) It is not necessary that $$(\eta_{\i }(\cd), \z_{\i }(\cd),\z_{\i }^M(\cd))\in L_{\dbF}^{2}(0,\i;\dbR^n)\times L_{\dbF}^{2}(0,\i;\dbR^n)\times M_{\dbF_-}^{2}(0,\i;\dbR^n).$$

(2) It is not necessary that $\lim_{t\rightarrow\i}\dbE|\eta_{\sc\i}(t)|^2=0$. Therefore, the terminal condition disappears in  \eqref{BSDEIN}.
\end{remark}

With the help of $\eta_{\sc\i}(\cd)$, now we can define the following close-loop control
\begin{equation}\label{optin2}\bar u_{\sc\i}^{x,\imath}(t)=\Th_{\sc\i}(\a(t))X(t)+ v_{\sc\i}(t).\end{equation}
where
$$\begin{cases}&\Th_{\sc\i}(\imath)=-\wt R(t,\imath)^{-1}(B^\top(\imath) P_{\sc\i}(\imath)+ D^\top(\imath) P_{\sc\i}(\imath)C(\imath)+ S(\imath)),
\\
& v_{\sc\i}(t,\imath)=-\wt R(t,\imath)^{-1}(D^\top(\a(t)) P_{\sc\i}(\imath)\si(t)+B^\top(\imath) \eta_{\sc\i}(t)+D^\top(\imath) \z_{\sc\i}(t)+r(t)).\end{cases}$$
Then the corresponding state process by $\bar X^{x,\imath}_{\sc\i}(\cd)$ satisfies
\bel{closed-loopin}\left\{\2n\ba{ll}
\ds d\bar X_{\sc\i}^{x,\imath}(t)=[A^{\Th_{\sc\i}}(t,\a(t))\bar X_{\sc\i}^{x,\imath}(t)+B(\a(t))v_{\sc\i}(t,\a(t))+b(t)]dt\\
\ns\ds\qq\qq+[C^{\Th_{\sc\i}}(t,\a(t))\bar X_{\sc\i}^{x,\imath}(t)+D(\a(t))v_{\sc\i}(t,\a(t))+\si(t)]dW(t),\\
\ns\ds\bar X_{\sc\i}^{x,\imath}(0)=x,\qq\a(0)=\imath.\ea\right.\ee

Our key result lies in deriving the estimate between  $\bar X_{\sc\i}^{x,\imath}(\cd)$ and $\bar X_{\scT}^{x,\imath}(\cd)$ (see \eqref{closed-loop}). This will be carefully studied in the next section. Before the end of this section,  some estimates are presented in the following proposition. The proof is posted in Section \ref{sec:pro}.

\begin{proposition}\label{festimateprop}Let {\rm(H1)--(H4)} and {\rm(H6)$'$} hold. Then for any $t\in[0,T]$, we have
\begin{align}
&\dbE\int_0^t e^{-{\d\over4}(t-s)}[|\z_{\scT}(s)|^2+|\z_{\sc\i}(s)|^2]dt\les K\int_0^\i\3n e^{-{\d\over4}|t-s|}\xi(s)ds,\label{festimate1}\\
&\dbE\int_0^t e^{-{\d\over4}(t-s)}|\z_{\scT}(s)-\z_{\sc\i}(s)|^2ds\les K e^{-\frac\d 8(T-t)}\int_0^\i  e^{-{\d\over4}|t-s|}\xi(s)ds,\label{festimate2}\\
&\dbE\int_0^te^{-\frac\d4(t-s)}| v_{\scT}(s)-v_{\sc\i}(s)|^2 ds\les Ke^{-{\d\over8}(T-t)}\int_0^\i e^{-{\d\over4}|t-s|}\xi(s)dr\label{v-v},\\
&\dbE\int_0^T|\zeta_{\scT}(s)|^2+|\zeta_{\sc\i}(s)|^2dt\les K\int_0^T
\xi(s)ds+K (T+1)\sup_{s\ges 0}\int_{0}^\infty e^{-\frac\d 4|s-r|}\xi(r)dr,\label{XXb}
\\
&\dbE\int_0^T|\zeta_{\scT}(s)-\zeta_{\sc\i}(s)|^2ds\les K\sup_{s\ges 0}\int_{0}^\infty e^{-\frac\d 4|s-r|}\xi(r)dr,\label{XXb2}\\
&\dbE\big[|\bar X^{x,\imath}_{\scT}(t)|^2+|\bar X^{x,\imath}_{\sc\i}(t)|^2\big]\les K \(e^{-{\d\over2}t}|x|^2+\int_0^\i e^{-{\d\over4}|t-s|}\xi(s)ds\).\label{boundEXTX}
\end{align}
\end{proposition}
\section{Strong Turnpike Property}\label{sec:STP}

In this section, we are going to state and prove the main result of this paper.

\bt{thmcase1} Let {\rm (H1)--(H5)} and {\rm(H6)$'$} hold. Let $(\bar X_{\scT}^{x,\imath}(\cd),\bar u^{x,\imath}_{\scT }(\cd))$ be the open-loop optimal pair of Problem {\rm(LQ)$_{\scT}$} corresponding to $(x,\imath)\in\sD$ (see \eqref{closed-loop}), and let $(\bar X^{x_{\sc\i},\imath}_{\sc\i}(\cd),\bar u^{x_\i,\imath}_{\sc\i}(\cd))$ be the state-control pair corresponding initial couple $(x_{\sc\i},\imath)\in\sD$ so that $\bar u^{x_\i,\imath}_{\sc\i}(\cd)$ given by \eqref{optin2} (see \eqref{closed-loopin}). Then it follows that
\bel{X-X}\ba{ll}
\ns\ds\dbE(|\bar X_{\scT}^{x,\imath}(t)-\bar X_{\sc\i}^{x_\i,\imath}(t)|^2+\int_0^te^{-{\d\over4}(t-s)}|\bar u_{\scT}^{x,\imath}(s)-\bar u_{\sc\i}^{x_\i,\imath}(s)|^2ds\\
\ns\ds\les Ke^{-{\d\over4}t}|x_{\sc\i}-x|^2+Ke^{-{\d\over8}(T-t)}\(e^{-{\d\over4}t}|x|^2+\int_0^\i e^{-{\d\over4}|t-s|}\xi(s)ds\).\ea\ee
for all $t\in[0,T]$.

\et

Before presenting the proof, some observations should be made here.
Taking $t=0$ and $t=T$, the right-hand side of \rf{X-X} respectively reads
$$K\[|x_{\sc\i}-x|^2+e^{-{\d\over2}T}
\(|x|^2+\int_0^\i e^{-{\d\over2}r}\xi(r)dr\)\],$$
and
$$K\[e^{-{\d\over2}T}|x_{\sc\i}-x|^2+\(e^{-{\d\over2}T}|x|^2+\int_0^\i e^{-{\d\over2}|T-r|}\xi(r)dr\)\],$$
which might not be small. However, for any $\e\in(0,{1\over2})$, if $t\in[\e T,(1-\e)T]$ (a middle range of $[0,T]$), then the right-hand side of \rf{Si(X-X)} can be estimated as follows:
$$\ba{ll}
\ns\ds K\[e^{-{\d\over2}t}|x_{\sc\i}-x|^2+e^{-{\d\over2}(T-t)}
\(e^{-{\d\over2}t}|x|^2+\int_t^\i e^{-{\d\over2}r}\xi(r)dr\)\]\\
\ns\ds\les K\[e^{-{\d\over2}\e T}|x_{\sc\i}-x|^2+e^{-{\d\over2}\e T}
\(e^{-{\d\over2}\e T}|x|^2+\int_t^\i e^{-{\d\over2}r}\xi(r)dr\)\]\to0,\ea$$
as $T\to\i$. This exactly describes what we call the turnpike property of our LQ problem. Now let us turn to the proof.


\ms

\begin{proof}[Proof of Theorem \ref{thmcase1}] In what follows, we will suppress the superscript $(x,\imath)$ and $(x_{\sc\i},\imath)$, together with $(t,\a(t))$. Set
$$\left\{\2n\ba{ll}
\ns\ds\wh X(t)=\bar X_{\sc\i}(t)-\bar X_{\scT}(t),\q\wh x=x_{\sc\i}-x,\q\wh u(t)=\bar u_{\sc\i}(t)-\bar u_{\scT}(t),\\
\ns\ds\wh\Th(t)=\Th_{\sc\i}(\a(t))-\Th_{\scT}(t,\a(t)),\qq\wh v(t)=v_{\sc\i}(t,\a(t))-v_{\scT}(t,\a(t)),\ea\right.$$
with $(\bar X_{\scT}(\cd),\bar u_{\scT}(\cd))$ being defined in \eqref{closed-loop} and $(\bar X_{\sc\i}(\cd),\bar u_{\sc\i}(\cd))$ being defined in \eqref{closed-loopin}. Since
$$\wh u(t)=\Th_{\sc\i}\wh X(t)+\wh\Th(t)\bar X_{\scT}(t)+\wh v(t),$$
we need only to estimate $\wh X(\cd),\wh v(\cd)$ in certain sense and $\bar X_{\scT}(\cd)$ uniformly bounded. The rest of this paper aims to realize this. Note that
$$\ba{ll}
\ns\ds d\widehat X(t)=d[\bar X_{\sc\i}(t)-\bar X_{\scT}(t)]\\
\ns\ds=\[\(A^{\Th_{\sc\i}}\bar X_{\sc\i}(t)+Bv_{\sc\i}(t)+b(t)\)-\(A^{\Th_{\scT}}\bar X_{\scT}(t)+Bv_{\scT}(t)+b(t)\)\]ds\\
\ns\ds\q+\[\(C^{\Th_{\sc\i}}\bar X_{\sc\i}(t)+Dv_{\sc\i}(t)+\si(t)\)-\(C^{\Th_{\scT}}\bar X_{\scT}(t)+Dv_{\scT}(t)+\si(t)\)\]dW(t)\\
\ns\ds=\(A^{\Th_\i}\wh X(t)+B\wh\Th(t)
\bar X_{\scT}(t)+B\wh v(t)\)dt+\(C^{\Th_\i}\wh X(t)+D\wh\Th(t)\bar X_{\scT}(t)+D\wh v(t)\)dW(t).\ea$$
Now, let $\Si_{\sc\i}(\cd)\in\hTh$ satisfy \rf{Ly1}. Applying It\^o's formula to $t\mapsto\lan\Si_{\sc\i}(\a(t))\wh X(t),\wh X(t)\ran$, we have
$$\ba{ll}
\ns\ds{d\over ds}\dbE\lan\Si_{\sc\i}(\a(t))\wh X(t),\wh X(t)\ran\\
\ns\ds=\1n\dbE\lan\(\L[\Si_{\sc\i}]\1n+\1n\Si_{\sc\i}A^{\Th_{\i}}\1n
+\1n(A^{\Th_{\i}})^\top\Si_{\sc\i}\1n+\1n(C^{\Th_{\i}})^\top
\Si_{\sc\i}C^{\Th_{\i}}\)\wh X(t),\wh X(t)\ran\\
\ns\ds\qq+2\dbE\lan\Si_{\sc\i}[B\wh\Th(t)
\bar X_{\scT}(t)+B\wh v(t)],\wh X(t)\ran+2\dbE\lan\Si_{\sc\i}C^{\Th_{\sc\i}}\wh X(t),D\wh v(t)+D\wh\Th(t)\bar X_{\scT}(t)\ran
\\
\ns\ds\qq+\dbE\lan\Si_{\sc\i}[D\wh v(t)+D\wh\Th(t)\bar X_{\scT}(t)],D\wh v(t)+D\wh\Th(t)\bar X_{\scT}(t)\ran\\
\ns\ds\les-{\d\over2}\dbE\lan\Si_{\sc\i}(\a(t))\wh X(t),\wh X(t)\ran+K\dbE|\wh\Th(t)\bar X_{\scT}(t)|^2+K\dbE|\wh v(t)|^2.\ea$$
By Gronwall's inequality, using \eqref{boundEXTX} and \eqref{v-v},
we have
\bel{Si(X-X)}\ba{ll}
\ns\ds\dbE|\wh X(t)|^2\les K\dbE\lan\Si_{\sc\i}(\a(t))\wh X(t),\wh X(t)\ran\\
\ns\ds\les Ke^{-{\d\over2}t}\lan\Si(\imath)\wh x,\wh x\ran+
K\int_0^te^{-{\d\over2}(t-s)}\dbE\(|\wh\Th(s)\bar X_{\scT}(s)|^2+|\wh v(s)|^2\)ds\\
\ns\ds\les Ke^{-{\d\over2}t}|\wh x|^2+K\int_0^te^{- {\d\over2}(t-s)}\dbE\(e^{-{\d\over2}(T-s)}|\bar X_{\scT}(s)|^2+|\wh v(s)|^2\)ds\\
\ns\ds\les Ke^{-{\d\over2}t}|\wh x|^2+K\int_0^te^{- {\d\over2}(t-s)}\dbE|\wh v(s)|^2ds\\
\ns\ds\qq+K\int_0^te^{- {\d\over2}(t-s)}e^{-{\d\over2}(T-s)}\(e^{-{\d\over2}s}|x|^2+\int_0^\i e^{-{\d\over4}|s-r|}\xi(r)dr\)ds\\
\ns\ds\les  Ke^{-{\d\over2}t}|\wh x|^2+Ke^{-{\d\over8}(T-t)}\(e^{-{\d\over4}t}|x|^2+\int_0^\i e^{-{\d\over4}|t-s|}\xi(s)ds\).\ea\ee
Note that
\bel{redU}\ba{ll}
\ns\ds\dbE|\wh u(s)|=\dbE|\Th_{\scT}(s)\bar X_{\scT}(s)-\Th_{\sc\i}(s)\bar X_{\sc\i}(s)+\bar v_{\scT}(s)-\bar v_{\sc\i}(s)|\\
\ns\ds\les\dbE\(|\Th_{\sc\i}(s)|\,|\wh X(s)|+|\wh\Th(s)|\,|\bar X_{\scT}(s)|+|\wh v(s)|\).\ea\ee
Using the estimates obtained in \eqref{Si(X-X)}, \eqref{boundEXTX}  and \eqref{v-v}, it follows that
\begin{align}&\int_0^te^{-{\d\over4}(t-s)}|\wh u(s)|^2ds=\int_0^te^{-{\d\over4}(t-s)}|\bar u_{\scT}^{x,\imath}(s)-\bar u_{\sc\i}^{x_\i,\imath}(s)|^2ds\nonumber\\
&\les K\int_0^te^{-{\d\over4}(t-s)}\dbE\(|\Th_{\sc\i}(s)|^2|\wh X(s)|^2+|\wh\Th(s)|^2|\bar X_{\scT}(s)|^2+|\wh v(s)|^2\)ds\nonumber\\
&\les K\int_0^te^{-{\d\over4}(t-s)}\dbE\(|\Th_{\sc\i}(s)|^2|\wh X(s)|^2+|\wh\Th(s)|^2|\bar X_{\scT}(s)|^2+|\wh v(s)|^2\)ds\nonumber\\
&\les Ke^{-{\d\over4}t}|\wh x|^2+Ke^{-{\d\over8}(T-t)}\(e^{-{\d\over4}t}|x|^2+\int_0^\i e^{-{\d\over4}|t-s|}\xi(s)ds\).\label{deltauu}
\end{align}
Our main result holds from \eqref{Si(X-X)} and \eqref{deltauu}.
\end{proof}

Until now, we have proven the STP for the optimal pair in Problem (LQ)$_{\scT}$. We can see the limit pair $(\bar X_{\sc\i}^{x,\imath}(\cd),\bar u_{\sc\i}^{x,\imath}(\cd))$ is identified by taking $T\rightarrow\infty$ for $(\Th_{\scT}(\cd), v_{\scT}(\cd))$.  Naturally, the next is to verify the optimality of the $(\bar X_{\sc\i}^{x,\imath}(\cd),\bar u_{\sc\i}^{x,\imath}(\cd))$ in some appropriate sense.

\section{Optimality of $(\bar X_{\sc\i}^{x,\imath}(\cd),\bar u_{\sc\i}^{x,\imath}(\cd))$}\label{sec:opt}
In this section, we will construct the appropriate optimal control problems for which $(\bar X_{\sc\i}^{x,\imath}(\cd),\bar u_{\sc\i}^{x,\imath}(\cd))$ is the optimal couple. We have two different cases.

\subsection{Integrable Case}\label{sec:int}

In this subsection, we work with integrable cases by assuming\ms

{\bf (H7).} $b(\cd),\si(\cd),q(\cd)\in L_{\dbF}^2(0,\i;\dbR^n),\q r(\cd)\in L_{\dbF}^2(0,\i;\dbR^m).$

\ms

It is obvious that (H7) is stronger than  (H6)$'$.   Recall that
the definition of $\xi(\cd)$ in \eqref{defxi}, $\xi(\cd)$ is integrable on $[0,\infty)$ with
$$\int_0^\i e^{-\beta|s-t|}\xi(s)ds\les\int_0^\i \xi(s)ds<\infty,$$
for any $\beta>0.$ All the previous results are true.


Recall that
\begin{equation*}\bar u_{\sc\i}^{x,\imath}(t)=\Th_{\sc\i}(\a(t))\bar X_{\sc\i}^{x,\imath}(t)+\bar v_{\sc\i}(t,\a(t)).\end{equation*}
where
$$\left\{\2n\ba{ll}
\ds\Th_{\sc\i}=-(R+D^\top P_{\sc\i}D)^{-1}(B^\top P_{\sc\i}+ D^\top P_{\sc\i}C+ S)\in\BS[A,C;B,D],\\
\ns\ds\bar v_{\sc\i}(t,\imath)=-(R+D^\top P_{\sc\i}D)^{-1}(D^\top P_{\sc\i}\si+B^\top \eta_{\sc\i}+D^\top \z_{\sc\i}+r).\ea\right.$$
We will see that \eqref{optin2} is the optimal control for a LQ problem on an infinite horizon.

To define the problem, we need the following set of admissible controls
$$\sU_{ad}^{x,\imath}[0,\i]=\Big\{ u(\cd)\in L_{\dbF}^2(0,\i;\dbR^m)\Big|  X(\cd;x,\imath,u(\cd))\in L_{\dbF}^2(0,\i;\dbR^n)
\Big\}$$
where $ X(\cd;x,i,u(\cd))$ is the solution of \eqref{state} with initial $(x,\imath)$ and control $u(\cd).$ For each $u(\cd)\in \sU_{ad}^{x,\imath}[0,\i]$, we define the following cost functional
$$J_{\sc\i}(x,\imath; u(\cd))=\dbE\(\int_0^\i g(t,
X(t),\a(t),u(t))dt\).$$
It can be easily seen that the cost functional is well-defined. We have the following LQ optimization problem on $[0,\i]$.\ms

{\bf Problem (LQ)$_{\i}$.} For a given initial $(x,\imath)\in\sD$, find a control $\bar u_{\scT}^{x,\imath}(\cd)\in\sU_{ad}^{x,\imath}[0,\i]$ such that
\bel{opt-u4}J_{\sc\i}(x,\imath;\bar u_{\scT}^{x,\imath}(\cd))=\inf_{u(\cd)\in\sU_{ad}^{x,\imath}
[0,\i]}J_{\i}(x,\imath;u(\cd))\equiv V_{\sc\i}(x,\imath).\ee

Now let us verify the optimality of \eqref{optin2} for Problem (LQ)$_{\i}$ in the following proposition.

\begin{proposition} Under {\rm (H1)-(H5)}, {\rm (H6)$'$} and {\rm (H7)}, $(\bar X_{\sc\i}^{x,\imath}(\cd),\bar u_{\sc\i}^{x,\imath}(\cd))$ is the unique optimal pair for Problem (LQ)$_{\i}$.
\end{proposition}
The above proposition is a special case studied in \cite{Mei-Wei-Yong-2025} and hence the proof is omitted. Now we can conclude the following corollary immediately.

\begin{corollary} Under the same assumptions as in Theorem \ref{thmcase1}, it follows that
\begin{align*}
\limsup_{T\rightarrow\infty}\(\dbE\int_0^T|\bar X_{\scT}^{x,\imath}(t)-\bar X_{\sc\i}^{x,\imath}(t)|^2 +|\bar u_{\scT}^{x,\imath}(t)-\bar u_{\sc\i}^{x,\imath}(t)|^2\)dt=0.
\end{align*}
\end{corollary}

\begin{proof} Write $h(t)=\int_0^{\infty}e^{-\frac{\d}4|s-t|}\xi(s)ds$. It is straightforward to see that (H7) yields that
$\int_0^\i h(t)dt<\infty.$
Note that
\begin{align*}
&\int_0^Th(t)e^{-\frac{\d}8(T-t)}dt=\int_0^{T/2}h(t)e^{-\frac{\d}8(T-t)}dt+\int_{T/2}^Th(t)e^{-\frac{\d}8(T-t)}dt\\
&\les e^{-\frac{\d}{16}T}\int_0^{\infty}h(t)dt+\int_{T/2}^\infty h(t)dt\rightarrow 0,\text{ as } T\rightarrow\infty.
\end{align*}
Then we have
$$\limsup_{T\rightarrow\infty}\dbE\int_0^T|\bar X_{\scT}^{x,\imath}(t)-\bar X_{\sc\i}^{x,\imath}(t)|^2dt=0.$$
By \eqref{redU}, it follows that
\begin{align*}&\limsup_{T\rightarrow\infty}\dbE\int_0^T|\bar u_{\scT}^{x,\imath}(t)-\bar u_{\sc\i}^{x,\imath}(t)|^2dt\\
&\les K \limsup_{T\rightarrow\infty}\dbE\int_0^T|\bar X_{\scT}^{x,\imath}(t)-\bar X_{\sc\i}^{x,\imath}(t)|^2+e^{\d(T-t)}\dbE|\bar X_{\scT}^{x,\imath}(t)|^2+|v_{\scT}(t)-v_{\sc\i}(t)|^2dt\\
&= K\limsup_{T\rightarrow\infty}\dbE\int_0^T|v_{\scT}(t)-v_{\sc\i}(t)|^2dt\\
&\les K\limsup_{T\rightarrow\infty}\dbE\int_0^T|\eta_{\scT}(t)-\eta_{\sc\i}(t)|^2+|\z_{\scT}(t)-\z_{\sc\i}(t)|^2+e^{-\d(T-t)}\(|\zeta_{\scT}(t)|^2+|\eta_{\scT}(t)|^2+\xi(t)\)dt\\
&=0.\end{align*}
In the above, we have used the boundedness of $\dbE|\bar X_{\scT}^{x,\imath}(t)|^2$ and $\dbE|\bar \eta_{\scT}(t)|^2$ (see \eqref{boundEXTX} and \eqref{stabsde}), \eqref{E[eta-eta]}, \eqref{festimate1}, \eqref{XXb2} and (H6)$'$.
\end{proof}

Before we finish this subsection, it is worth remarking that even without the switching Markov chain, our results are not studied in  \cite{Sun-Wang-Yong-2022} or \cite{Sun-Yong-2024S} where $b,\sigma,q,r$ are assumed to be deterministic and constants. Our assumption (H7) allows those non-homogeneous terms to be stochastic. With some appropriate integrability conditions, we derive a new form of STP compared to that in \cite{Sun-Wang-Yong-2022} and  \cite{Sun-Yong-2024S}.

\subsection{Local-Integrable Case}\label{sec:loc}
In this subsection, we work with local-integrable cases by assuming
\ms

{\bf (H8)} $b(\cd),\si(\cd),q(\cd)\in L_{\dbF}^{2}(0,T;\dbR^n),\q r(\cd)\in L_{\dbF}^{2}(0,T;\dbR^m)$ for each $T>0$ with
\begin{equation}\label{boundzeta}\limsup_{T\rightarrow\infty}\frac1T\int_0^T\xi(t)dt<\infty.\end{equation}

In this section, we always assume that (H1)-(H5), (H6)$'$ and (H8) hold.  In this case, we can show that $\bar u_{\sc\infty}^{x,\imath}(\cd)$ through the limit process is the optimal control for the following ergodic control problem.
Define $$\sU_{loc}[0,\i)=\bigcap_{T\geq 0}\sU[0,T].$$

{\bf Problem (LQ)$_E$.}
For a given initial state $(x,\imath)\in\sD$, find a control $u_E(\cd)\in\sU_{loc}[0,\i)$ such that
\bel{opt-u2}
J_{\rm E}(x,\imath;u_E(\cd))=\inf_{u(\cd)\in\sU_{loc}[0,\i)}J_E(x,\imath;u(\cd))=: V_{\rm E}(x,\imath),\ee
where the ergodic cost is defined by
\begin{equation*}
J_E(x,\imath;u(t)):=\liminf_{T\rightarrow\infty}\frac1TJ_{\rm E}(x,\imath;u_E(\cd)).
\end{equation*}

\begin{proposition} Suppose  {\rm (H1)--(H5)}, {\rm (H6)$'$} and {\rm (H8)} hold. For any $(x,\imath)\in\sD$, $\bar    u^{x,\imath}_{\sc\infty}(\cd)$ is the optimal control and $\bar X^{x,\imath}_{\sc\infty}(\cd)$ is the corresponding optimal trajectory   for Problem (LQ)$_E$.
Moreover, $J_E(x,\imath;\bar u^{x,\imath}_{\sc\infty}(\cd))$ is finite.
\end{proposition}

\begin{proof} We suppress the top index $(x,\imath)$ in the proof. Note that
\begin{align*}&|\bar u_{\scT}(t)|+|\bar u_{\i}(t)|\\
&\les K\(|X_{\scT}(t)|+|\eta_{\scT}(t)|+|\zeta_{\scT}(t)|+|X_{\sc\i}(t)|+|\eta_{\sc\i}(t)|+|\zeta_{\sc\i}(t)|+|r(t)|+|\sigma(t)|\)\end{align*}
and
\begin{align*}&|\bar u_{\scT}(t)-\bar u_{\i}(t)|\les |\Th_{\sc\i}(\bar X_{\scT}(t)-X_{\sc\i}(t))|+|(\Th_{\sc\i}-\Th_{\scT})X_{\scT}(t))|+|v_{\scT}(t)-v_{\sc\i}(t)|\\
&\les K|\bar X_{\scT}(t)-X_{\sc\i}(t)|+K|\eta_{\scT}(t)-\eta_{\sc\i}(t)|+K|\zeta_{\scT}(t)-\zeta_{\sc\i}(t)|\\
&\q+Ke^{-\frac\d2(T-t)}\(|X_{\scT}(t)|+|\eta_{\scT}(t)|+|\zeta_{\scT}(t)|+|r(t)|+|\sigma(t)|\)\end{align*}
Applying all the estimates in Proposition \ref{festimateprop}, \eqref{E[eta-eta]} and \eqref{X-X}, it follows that
\begin{align}\label{bdinT}&\limsup_{T\rightarrow\infty}\frac1T\dbE\int_0^T|\bar u_{\scT}(t)|^2+|\bar u_{\i}(t)|^2dt<\i\text{ and }\lim_{T\rightarrow\infty}\frac1T\dbE\int_0^T|\bar u_{\scT}(t)-\bar u_{\i}(t)|^2dt=0.\end{align}

 Next, we see
\begin{align*}
&\frac1T\dbE\int_0^ Tg(t,\a(t),X(t),u(t))dt\ges \frac1T\dbE\int_0^ Tg(t,\a(t),\bar X_{\scT}(t),\bar u_{\scT}(t))dt\\
&\ges\frac1T\dbE\int_0^ Tg(t,\a(t),\bar X_{\sc\i}(t),\bar u_{\sc\i}(t))dt\\
&\q-\frac KT \int_0^ T\Big(\dbE[|\bar X_{\scT}(t)-\bar X_{\sc\i}(t)|^2+|\bar u_{\scT}(t)-\bar u_{\sc\i}(t)|^2]\\
&\qq\qq\qq\cdot\dbE[1+|\bar X_{\scT}(t)|^2+|\bar X_{\sc\i}(t)|^2+|\bar u_{\scT}(t)|^2+|\bar u^{x,\imath}_{\sc\i}(t)|^2]\Big)^{\frac12}dt\\
&\q-\frac KT \int_0^ T\Big(\dbE[|\bar X_{\scT}(t)-\bar X_{\sc\i}(t)|^2+|\bar u_{\scT}(t)-\bar u_{\sc\i}(t)|^2]\Big)^{\frac12}\\
&\qq\qq\qq\cdot\(\dbE[1+|\bar X_{\scT}(t)|^2+|\bar X_{\sc\i}(t)|^2+| \bar u_{\sc\infty}(t)|^2+| \bar u_{\scT}(t)|^2]\)^{\frac12}dt.
\end{align*}
Taking $T\rightarrow\infty$, it follows that for any $u(\cd)\in\sU_{loc}[0,\infty)$,
\begin{align*}&\liminf_{T\rightarrow\infty}\frac1T\dbE\int_0^ Tg(t,\a(t),X(t),u(t))dt\ges\liminf_{T\rightarrow\infty}\frac1T\dbE\int_0^ Tg(t,\a(t),\bar X^{x,\imath}_{\scT}(t),\bar u^{x,\imath}_{\scT}(t))dt\\
&=\liminf_{T\rightarrow\infty}\frac1T\dbE\int_0^ Tg(t,\a(t),\bar X^{x,\imath}_{\sc\i}(t),\bar u^{x,\imath}_{\sc\i}(t))dt. \end{align*}
Moreover, the uniform boundedness of $\dbE|\bar X_{\sc\i}|^2$ and \eqref{bdinT} together imply that
\begin{align*}&\frac1 T\int_0^T\dbE|\bar X_{\sc\i}(t)|^2dt<K,\q\frac1 T\int_0^T\dbE|\bar u_{\sc\i}(t)|^2dt<K,\\
&\frac1 T\int_0^T\dbE|\lan q(t),\bar X_{\sc\i}(t)\ran|dt\les \frac1T\int_0^T\xi(t)+\dbE|\bar X_{\sc\i}(t)|^2dt<K,\\
&\frac1 T\int_0^T\dbE|\lan r(t),\bar u_{\sc\i}(t)\ran|dt\les \frac1T\int_0^T\xi(t)+\dbE|\bar u_{\sc\i}(t)|^2dt<K.
\end{align*}
Therefore $J_E(x,\imath;\bar u_{\sc\i}(t))$ is finite. Moreover, $\bar u_{\sc\i}(t)$ is the optimal control process in $\sU_{loc}[0,\i)$ and $\bar X_{\sc\i}(t)$ is the corresponding  trajectory for Problem (LQ)$_E$.
\end{proof}



Finally, when $b(\cd),\sigma(\cd),q(\cd),r(\cd)$ are bounded and $\dbF^\a$-measurable (instead of $\dbF$-measurable), one can easily see that $\z_{\scT}(t)=0$ for all $T>0$ and $0\les t\les T$. Then all the estimates on $\zeta_{\scT}$ are not necessary. For such a particular case, without essential difficulties, \eqref {X-X} can be refined to \bel{X-X2}\ba{ll}
\ns\ds\dbE(|\bar X_{\scT}^{x,\imath}(t)-\bar X_{\sc\i}^{x_\i,\imath}(t)|^2+|\bar u_{\scT}^{x,\imath}(t)-\bar u_{\sc\i}^{x_\i,\imath}(t)|^2\\
\ns\ds\les Ke^{-{\d\over2}t}|x_{\sc\i}-x|^2+Ke^{-{\d\over8}(T-t)}\(e^{-{\d\over4}t}|x|^2+\int_0^\i e^{-{\d\over4}|t-s|}\xi(s)ds\).\ea\ee
The above matches  the results obtained in \cite{Sun-Wang-Yong-2022} and \cite{Sun-Yong-2024S} where $b(\cd),\sigma(\cd), q(\cd),r(\cd)$ are assumed to be deterministic constants. From this, we can see that those previous results in \cite{Sun-Wang-Yong-2022} and \cite{Sun-Yong-2024S} are some special cases of those in the current paper, even without switching states.

\section{Concluding Remarks}\label{sec:con}
In this paper, we obtained the turnpike property for LQ optimal control in an infinite horizon with a regime-switching state when the system is non-homogeneous. We see that the comparing limit pair admits different optimality for different optimal control problems,  depending on the integrability of the optimal solution over the infinite horizon. Those relate to three different cases: homogeneous case, integrable case, and local-integrable case. Even for the problem without switching, our results provide more accurate bounds under weaker assumptions compared to the previous results in the literature.

\section{Proofs}\label{sec:pro}

In this section, we present the proofs of some results.

\begin{proof} [Proof of Proposition \ref{3.5}]




Due to the linearity of the BSDE, the existence and uniqueness of the adapted solution triple $(\eta_{\scT}(\cd),\z_{\scT}(\cd),\z_{\scT}^M(\cd))\in L_\dbF^2(0,T;\dbR^n)\times L_\dbF^2(0,T;\dbR^n)\times M_{\dbF_-}^2(0,T;\dbR^n)$ is standard. Thus, we only need to establish the estimates. We split the proof into several steps.

\ms

\it Step 1. Dissipativity of the modified system. \rm Note that the (homogenous) closed-loop system $[A^{\Th_{\scT}},C^{\Th_{\scT}}]$ may not be dissipative for some states of the Markov chain (see \rf{closed-loop}). 
The key step in our proof is to seek a modified and equivalent system that is dissipative for all the states of the Markov chain. Then the estimates can be derived in a classical way. To this end, let $\Si_{\sc\i}(\imath)\in\hSi$ be that in Proposition \ref{3.4} (i), satisfying \rf{Ly1}. Set
$$E(\imath,\jmath)=\Si_{\sc\i}(\jmath)^{1\over2}-\Si_{\sc\i}
(\imath)^{1\over2}\in\dbS^n.$$
It is clear that $E(\imath,\jmath)$ is well-defined and symmetric. By \rf{dF}, we have
\bel{dF*}d[\Si_{\sc\i}(\a(t))^{1\over2}]=\L[\Si_{\sc\i}^{1\over2}](\a(t))dt
+\sum_{\imath,\jmath\in\cM}E(\imath,\jmath){\bf1}_{\{\a(t^-)=\imath\}}dM_{\imath\jmath}(t),\ee
Let $X(\cd)$ be the solution of the homogeneous system $[A^{\Th_T},C^{\Th_T}]$. We define
\bel{wtX}\wt X_{\scT}(t)=\Si_{\sc\i}(\a(t))^{1\over2}X_{\scT}(t),\qq t\in[0,T].\ee
Then, by It\^o's formula, we have
$$\ba{ll}
\ns\ds d\wt X_{\scT}(t)=d[\Si_{\sc\i}(\a(t))^{1\over2}X_{\scT}(t)]\\
\ns\ds=\(\L[\Si_{\sc\i}^{1\over2}](\a(t))X_{\scT}(t)+\Si_{\sc\i}
(\a(t))^{1\over2}A^{\Th_T}(t,\a(t))X_{\scT}(t)\)dt\\
\ns\ds\qq+\Si_{\sc\i}(\a(t))^{1\over2}C^{\Th_{\scT}}(t,\a(t))X_{\scT}(t)dW(t)
+\sum_{\imath,\jmath\in\cM}E(\imath,\jmath)X_{\scT}(t^-)
{\bf1}_{\{\a(t^-)=\imath\}}dM_{\imath\jmath}(t)\\
\ns\ds\equiv\(\L[\Si_{\sc\i}^{1\over2}](\a(t))\Si_{\sc\i}(\a(t))^{-{1\over2}}
+\Si_{\sc\i}(\a(t))^{1\over2}A^{\Th_T}(t,\a(t))\Si_{\sc\i}(\a(t))
^{-{1\over2}}\)\wt X_{\scT}(t)dt\\
\ns\ds\qq+\Si_{\sc\i}(\a(t))^{1\over2}C^{\Th_{\scT}}(t,\a(t))\Si_{\sc\i}(\a(t))^{
-{1\over2}}\wt X_{\scT}(t)dW(t)\\
\ns\ds\qq+\sum_{\imath,\jmath\in\cM}E(\imath,\jmath)\Si_{\sc\i}(\a(t))^{
-{1\over2}}\wt X_{\scT}(t^-){\bf1}_{\{\a(t^-)=\imath\}}dM_{\imath\jmath}(t)\\
\ns\ds=\wt A^{\Th_T}\wt X(t)dt+\wt C^{\Th_T}\wt X(t)dW(t)+\sum_{\imath,\jmath\in\cM}\wt E(\imath,\jmath)\wt X_{\scT}(t^-){\bf1}_{\{\a(t^-)=\imath\}}dM_{\imath\jmath}(t),\ea$$
where
\bel{wtA}\ba{ll}
\ns\ds\wt A^{\Th_T}(t,\imath)=\L[\Si_{\sc\i}^{1\over2}](\imath)
\Si_{\sc\i}(\imath)^{-{1\over2}}+\Si_{\sc\i}(\imath)^{1\over2}A^{\Th_T}(t,\imath)
\Si_{\sc\i}(\imath)^{-{1\over2}},\\
\ns\ds\wt C^{\Th_T}(t,\imath)=\Si_{\sc\i}(\imath)^{1\over2}C^{\Th_T}
(t,\imath)\Si_{\sc\i}(\imath)^{-{1\over2}},\q\wt E(\imath,\jmath)=E(\imath,\jmath)\Si_{\sc\i}(\imath)^{-{1\over2}}.\ea\ee
Thus, we obtain the following new SDE
\bel{}\left\{\2n\ba{ll}
\ds d\wt X_{\scT}(t)\1n=\2n\wt A^{\Th_T}\1n\wt X(t)dt\1n+\1n\wt C^{\Th_T}\1n\wt X(t)dW(t)\1n+\3n\sum_{\imath,\jmath\in\cM}\3n\wt E(\imath,\jmath)\wt X_{\scT}(t^-){\bf1}_{(\a(t^-)=\imath)}
dM_{\imath\jmath}(t),\\
\ns\ds\qq\qq\qq\qq\qq\qq\qq\qq\qq\qq\qq t\in[0,T],\\
\ns\ds\wt X_{\scT}(0)=\Si(\imath)^{1\over2}x.\ea\right.\ee
Note that $|\wt X_{\scT}(t)|^2=\lan\Si_{\sc\i}(\a(t))X_{\scT}(t),X_{\scT}(t)\ran$, and by It\^o's formula, we have
$$\ba{ll}
\ns\ds{d\over dt}\dbE|\wt X(t)|^2=\dbE\Big\lan\Big(\wt A^{\Th_T}+(\wt A^{\Th_T})^\top+(\wt C^{\Th_T})^\top\wt C^{\Th_T}\\
\ns\ds\qq\qq\qq\qq\qq+\sum_{\jmath\ne\a(t)}\l_{\a(t)\jmath}\wt E(\a(t),\jmath)^\top\wt E(\a(t),\jmath)\Big)\wt X_{\scT}(t),\wt X_{\scT}(t)\Big\ran\\
\ns\ds={d\over dt}\lan\Si_{\sc\i}(\a(t))X_{\scT}(t),X_{\scT}(t)\ran\\
\ns\ds=\dbE\Big\lan\(\L[\Si_{\sc\i}](\a(t))+\Si_{\sc\i}(\a(t))
A^{\Th_{\scT}}+(A^{\Th_{\scT}})^\top\Si_{\sc\i}(\a(t))\\
\ns\ds\qq+(C^{\Th_{\scT}})^\top\Si_{\sc\i}(\a(t))C^{\Th_{\scT}}\)X_{\scT}(t),X_{\scT}(t)\Big\ran\\
\ns\ds\les-{\d\over2}\dbE\lan\Si_{\sc\i}(\a(t))X_{\scT}(t),X_{\scT}(t)\ran
=-{\d\over2}|\wt X_{\scT}(t)|^2,\qq t\in[0,T-T_0].\ea$$
Thus, we have for $t\in[0,T-T_0]$,
\bel{dissipative**}\wt A^{\Th_T}+(\wt A^{\Th_T})^\top+(\wt C^{\Th_T})^\top\wt C^{\Th_T}+\sum_{\jmath\ne\imath}\l_{\imath\jmath}\wt E(\imath,\jmath)^\top\wt E(\imath,\jmath)\les-{\d\over2}I.\ee
This means $\wt X_{\scT}(\cd)$ itself is dissipative on $[0,T-t_0]$, which will be very important below.

\ms

\it Step 2. BSDEs for the modified adapted solution. \rm Again, let $\Si_{\sc\i}(\cd)\in\hSi$ be that in Proposition \ref{3.4} (i), satisfying \rf{Ly1}. Let $(\eta_{\scT}(\cd),\z_{\scT}(\cd),\z^M_{\scT}(\cd))$ be the adapted solution to BDSE \rf{BSDET}. 
Now, we write
\bel{wteta}\ba{ll}
\ns\ds\wt\eta_{\scT}(t)=\Si_{\sc\i}(\a(t))^{-{1\over2}}\eta_{\scT}(t),\q \wt\z_{\scT}(t)=\Si_{\sc\i}(\a(t))^{-{1\over2}}\z_{\scT}(t),\\
\ns\ds\wt\z^M_{\scT}(t)=\Si_{\sc\i}(\a(t))^{-{1\over2}}\z^M_{\scT}(t),\q
\wt\f_{\scT}(t)=\Si_{\sc\i}(\a(t))^{-{1\over2}}\f_{\scT}(t,\a(t)),\ea
\q t\in[0,T],\ee
Then we claim that $(\wt\eta_{\scT}(\cd),\wt\z_{\scT}(\cd),\wt\z_{\scT}^M(\cd))$ is the adapted solution of the following BSDE (compared with \rf{BSDET}):
\bel{BSDEwteta}\left\{\2n\ba{ll}
\ds d\wt\eta_{\scT}(t)=-\(\wt A^{\Th_T}(t,\a(t))^\top\wt\eta_{\scT}(t)+\wt C^{\Th_T}(t,\a(t))^\top\wt\z_{\scT}(t)+\wt\f_{\scT}(t)\\
\ns\ds\qq\qq+\sum_{\imath,\jmath\in\cM} \wt E(\imath,\jmath)^\top\wt\z^M_{\scT}(t,\jmath)
\l_{\imath\jmath}{\bf1}_{\{\a(t)=\imath\}}\)dt+\wt\z_{\scT}(t) dW(t)+\wt\z^M_{\scT}(t)dM(t),\\
\ns\ds\qq\qq\qq\qq\qq\qq\qq\qq\qq\qq\qq t\in[0,T],\\
\ns\ds\wt\eta_{\scT}(T)=\vartheta.\ea\right.\ee
Here $\vartheta$ is a $\cF_T$ measurable random variable with finite second moment.In fact, noting $\eta_{\scT}(t)=\Si_{\sc\i}(\a(t))^{1\over2}\wt\eta_{\scT}(t)$, using It\^o's formula, we have
\begin{align*} &d\eta_{\scT}(t)=\(d[\Si_{\sc\i}(\a(t))^{1\over2}]\)\wt \eta_{\scT}(t^-)\1n+\1n\Si_{\sc\i}(\a(t^-))^{1\over2}d\wt\eta_{\scT}(t)\\
&\q+\3n\sum_{\imath,\jmath\in\cM}\3n E(\imath,\jmath)\wt \z_{\scT}^M(t,\jmath)\lambda_{\imath\jmath}{\bf1}_{\{\a(t)=\imath\}}dt
+\sum_{\imath,\jmath\in\cM}E(\imath,\jmath)\wt \z_{\scT}^M(t,\jmath){\bf1}_{\{\a(t^-)=\imath\}}dM_{\imath\jmath}(t)\\
&=\Big(\L[\Si_{\sc\i}(\cd)^{1\over2}](\a(t))\wt\eta_{\scT}(t) dt+\sum_{\imath,\jmath\in\cM}E(\imath,\jmath){\bf1}_{\{\a( t^-)=\imath\}}\wt
\eta(t^-)dM_{\imath\jmath}(t)\Big)\\
&\q+\Si_{\sc\i}(\a(t))^{1\over2}\[-\((\wt A^{\Th_T})^\top\wt\eta_{\scT}(t)\1n+\1n(\wt C^{\Th_T}_{\scT})^\top\wt\z_{\scT}(t)\1n+\1n\wt\f_{\scT}(t)\\
&\q+\3n\sum_{\imath,\jmath\in\cM}\3n \wt E(\imath,\jmath)^\top\wt\z^M_{\scT}(t,\jmath)\l_{\imath\jmath}
{\bf1}_{\{\a(t)=\imath\}}\1n\)dt\1n+\1n\wt \z_{\scT}(t)dW(t)\1n+\3n\1n\sum_{\imath,\jmath\in\cM}\2n\wt \z^M_{\scT}(t,\jmath){\bf1}_{\{\a(t^-)=\imath\}}dM_{\imath\jmath}(t)\]\\
&\q+\sum_{\imath,\jmath\in\cM} E(\imath,\jmath)\wt \z^M_{\scT}(t,\jmath)\l_{\imath\jmath}{\bf1}_{\{\a(t)=\imath\}}dt
+\sum_{\imath,\jmath\in\cM}E(\imath,\jmath)\wt \z^M_{\scT}(t,\jmath){\bf1}_{\{\a(t^-)=\imath\}}dM_{\imath\jmath}(t)\\
&=-\[-\L[\Si_{\sc\i}^{1\over2}](\a(t))\wt\eta_{\scT}(t)
\1n+\1n\Si(\a(t))^{1\over2}\((\wt A^{\Th_T})^\top\wt\eta_{\scT}(t)\1n+\1n(\wt C^{\Th_T})^\top\wt\z_{\scT}(t)\\
&\q+\wt\f_{\scT}(t)\)dt\]+\Si_{\sc\i}(\a(t))^{1\over2}
\wt\z_{\scT}(t)dW(t)\\
&\q+\3n\sum_{\imath,\jmath\in\cM}\3n\Big(\Si_{\sc\i}(\jmath)^{1\over2}
\wt\z^M_{\scT}(t,\jmath)\1n+\1n E(\imath,\jmath)\wt \eta_{\scT}(t^-)\){\bf1}_{(\a(t^-)=\imath)}dM_{\imath\jmath}(t)\\
&=-[(A^{\Th_T})^\top\eta_{\scT}(t)+(C^{\Th_T})^\top
\z_{\scT}(t)+\f_{\scT}(t,\a(t))]dt+\z_{\scT}(t)dW(t)+\z^M_{\scT}(t)dM(t),\end{align*}
where (note \rf{wteta})
$$\left\{\2n\ba{ll}
\ds\z_{\scT}(t)=\Si_{\sc\i}(\a(t))^{1\over2}\wt\z_{\scT}(t),\\
\ns\ds\z^M_{\scT}(t,\jmath)=\Si_\i(\jmath)^{1\over2}\wt\z^M_{\scT}
(t,\jmath)+E(\a(t^-),\jmath)\wt\eta_{\scT}(t^-).\ea\right.$$
In the above,
$$\ba{ll}
\ns\ds\Si_{\sc\i}(\a(t))^{-{1\over2}}\sum_{\imath,\jmath\in\cM}\wt E(\imath,\jmath)^\top{\bf1}_{\{\a(t)=\imath\}}=\sum_{\imath,\jmath\in\cM}
\Si_{\sc\i}(\imath)^{-{1\over2}}\(E(\imath,\jmath)\Si(\imath)^{1\over2}
\)^\top{\bf1}_{\{\a(t)=\imath\}},\\
\ns\ds=\sum_{\imath,\jmath\in\cM}E(\imath,\jmath)^\top{\bf1}_{\{\a(t)=\imath\}}
=\sum_{\imath,\jmath\in\cM}E(\imath,\jmath).\ea$$
By the uniqueness of a linear BSDE, the above calculation yields that $(\wt\eta_{\scT}(\cd),\wt\z_{\scT}(\cd),\wt\z_{\scT}^M(\cd))$ is the adapted solution of the BSDE \rf{BSDEwteta} by taking $\vartheta=0$.

\ms

\it Step 3. Dissipation inequality for $\wt\eta_{\scT}(\cd)$. \rm By It\^o's formula, for $t\in[0,T-T_0]$, we have
$$\ba{ll}
\ns\ds{d\over dt}\dbE|\wt\eta_{\scT}(t)|^2=\(-2\lan\wt\eta_{\scT}(t),
\wt A^{\Th_T}(t,\a(t))^\top\wt\eta_{\scT}(t)+\wt  C^{\Th_T}(t,\a(t))^\top\wt\z_{\scT}(t)+\wt\f_{\scT}(t)\ran\\
\ns\ds\qq\qq\qq+\sum_{\imath\ne\jmath} E(\imath,\jmath)\wt\z^M_{\scT}(t,\jmath)
\l_{\imath\jmath}{\bf1}_{\{\a(t)=\imath\}}\ran+|\wt\z_{\scT}(t)|^2dt
+\sum_{\imath\neq\jmath}\lambda_{\imath\jmath}|\wt\z^M_{\scT}(t,\jmath)|^2{\bf 1}_{[\a(t)=\imath]}\)dt\\
\ns\ds=-\dbE\(\lan[\wt A^{\Th_T}+(\wt A^{\Th_T})^\top]
\wt\eta_{\scT}(t),\wt\eta_{\scT}(t)\ran-|\wt C^{\Th_T}\wt\eta_{\scT}(t)|^2+|\wt\z_{\scT}(t)-\wt C^{\Th_T}(\a(t))\wt\eta_{\scT}(t)|^2\\
\ns\ds\qq+\3n\sum_{\imath\neq \jmath}\3n\(|\wt\z^M_{\scT}(t,\jmath)|^2
\2n-\1n2\lan E(\imath,\jmath)\wt\eta_{\scT}(t),\wt\z^M_{\scT}(t,\jmath)\ran\)
\l_{\imath\jmath}{\bf1}_{\{\a(t^-)=\imath\}}\2n-\1n2\dbE\lan\wt\f_{\scT}(t),
\wt\eta_{\scT}(t)\ran\)\\
\ns\ds=\1n-\dbE\[\1n\lan\(\wt A^{\Th_T}\2n+\1n(\wt A^{\Th_T})^\top\3n
+\1n(\wt C^{\Th_T})^\top\wt C^{\Th_T}\1n+\3n\sum_{\imath\neq \jmath}\2n E(\imath,\jmath)E(\imath,\jmath)
\l_{\imath\jmath}{\bf1}_{\{\a(t)=\imath\}}\)\wt\eta_{\scT}(t),
\wt\eta_{\scT}(t)\ran\\
\ns\ds\qq+|\wt\z_{\scT}(t)\1n-\1n\wt C^{\Th_T}\wt\eta_{\scT}(t)|^2
\2n+\2n\2n\sum_{\imath,\jmath\in\cM}\2n|\z^M_{\scT}(t,\jmath)
\1n-\1n E(\imath,\jmath)
\wt\eta_{\scT}(t)|^2
\l_{\imath\jmath}{\bf1}_{\{\a(t)=\imath\}}\2n-\1n2\dbE\lan\wt\f_{\scT}(t),
\wt\eta_{\scT}(t)\ran\1n\)\\
\ns\ds\ges\dbE\({\d\over4}|\wt\eta_{\scT}(t)|^2\2n+\2n|\wt\z_{\scT}(t)-\wt C^{\Th_T}(t,\a(t))\wt\eta_{\scT}(t)|^2\1n+\3n\sum_{\imath\neq \jmath}
|\wt\z^M_{\scT}(t,\jmath)\1n-\1n E(\imath,\jmath)\wt\eta_{\scT}(t)|^2
\l_{\imath\jmath}{\bf1}_{\{\a(t^-)=\imath\}}\\
\ns\ds\qq-K\dbE|\wt\f_{\scT}(t)|^2\).\ea$$
Where we have used the dissipativity of $\wt X(\cd)$, i.e., \rf{dissipative**}, in the last step. 
Thus, for all $t\in[0,T-T_0]$, the following dissipativity inequality holds:
\bel{deta}\ba{ll}
\ns\ds{d\over dt}\dbE|\wt\eta_{\scT}(t)|^2\ges\dbE\({\d\over4}|\wt\eta_{\scT}(t)|^2\2n+\2n|\wt\z_{\scT}(t)-\wt C^{\Th_T}(t,\a(t))\wt\eta_{\scT}(t)|^2\\
\ns\ds\qq\qq\qq\qq+\sum_{\imath\neq \jmath}
|\wt\z^M_{\scT}(t,\jmath)-E(\imath,\jmath)\wt\eta_{\scT}(s)|^2
\l_{\imath\jmath}{\bf1}_{\{\a(t)=\imath\}}-K\xi(t)\).\ea\ee

On $[T-T_0,T]$, using the boundedness of $A^{\Th_{\scT}}$ and $C^{\Th_{\scT}}$, it is standard to derive that
\bel{deta2}\ba{ll}
\ns\ds{d\over dt}\dbE|\wt\eta_{\scT}(t)|^2\ges \dbE\Big(-K|\wt\eta_{\scT}(t)|^2\2n+\2n|\wt\z_{\scT}(t)-\wt C^{\Th_T}(t,\a(t))\wt\eta_{\scT}(t)|^2\\
\ns\ds\qq\qq\qq\qq+\sum_{\imath\neq \jmath}
|\wt\z^M_{\scT}(t,\jmath)-E(\imath,\jmath)\wt\eta_{\scT}(s)|^2
\l_{\imath\jmath}{\bf1}_{\{\a(t)=\imath\}}-K\xi(t)\).\ea\ee

\it Step 4. Boundedness of $(\eta_{\scT}(\cd),\z_{\scT}(\cd),\z^M_{\scT}(\cd))$.  \rm By \eqref{deta2}, Grownwall's inequality implies that for $t\in[T-t_0, T],$
\begin{align*}
&\dbE|\wt \eta_{\scT}(t)|^2+\int_t^{T}e^{K(s-t)}\dbE|\wt \z_{\scT}(s)+ C^{\Th_{\scT}}(s,\a(s))\wt \eta_{\scT}(s)|^2ds\\
&\les K\int_t^{T}e^{K(s-t)}\xi(s)ds+e^{K(T-t)}\dbE|\vartheta|^2.
\end{align*}
Because $0\les s-t\les t_0$ for $t\in[T-t_0,T]$, the above is equivalent to
\begin{align}\label{wtYZ}
&\dbE|\wt \eta_{\scT}(t)|^2+\int_t^{T}e^{-\frac{\d}4(s-t)}\dbE|\wt \z_{\scT}(s)+C^{\Th_{\scT}}(s,\a(s))\wt \eta_{\scT}(s)|^2ds\nonumber\\
&\les K\int_t^{T}e^{-\frac{\d}4(s-t)}\xi(s)ds+e^{-\frac{\d}4(T-t)}\dbE|\vartheta|^2.
\end{align}
In particular, we have
\begin{align}\label{wtYZ-1}
&\dbE|\wt \eta_{\scT}(T-T_0)|^2+\int_{T-T_0}^{T}e^{-\frac{\d}4(s-T+T_0)}\dbE|\wt \z_{\scT}(s)+C^{\Th_{\scT}}(s,\a(s))\wt \eta_{\scT}(s)|^2ds\nonumber\\
&\les K\int_{T-T_0}^{T}e^{-\frac{\d}4(s-(T-T_0))}\xi(s)ds+K\dbE|\vartheta|^2.
\end{align}
%
By \eqref{deta} and \eqref{wtYZ-1}, Grownwall's inequality implies that for $t\in[0, T-T_0],$
\bel{Eeta}\ba{ll}
\ns\ds\dbE|\wt\eta_{\scT}(t)|^2+\int_t^{T-T_0}e^{-{\d\over4}(s-t)}\dbE|
\wt\z_{\scT}(s)-\wt C^{\Th_T}(s,\a(s))\wt\eta_{\scT}(s)|^2ds\\
\ns\ds+\int_t^{T-T_0}e^{-{\d\over4}(s-t)}\dbE\sum_{\imath,\jmath\in\cM}
|\wt\z^M_{\scT}(s,\jmath)-E(\imath,\jmath)\wt\eta_{\scT}(s)|^2
\l_{\imath\jmath}{\bf1}_{\{\a(s)=\imath\}}ds\\
\ns\ds\les e^{-{\d\over4}(T-T_0-t)}\dbE|\wt\eta_{\scT}(T-T_0)|^2+K\int_t^{T-T_0}
e^{-{\d\over4}(s-t)}\xi(s)ds\\
\ns\ds\les Ke^{-{\d\over4}(T-T_0-t)}\int_{T-T_0}^{T}e^{-\frac{\d}4(s-(T-T_0))}\zeta(s)ds+K\int_t^{T-T_0}
e^{-{\d\over4}(s-t)}\xi(s)ds+Ke^{-\frac{\d}4(T-t)}\dbE|\vartheta|^2\\
\ns\ds\les K\int_t^{T}
e^{-{\d\over4}(s-t)}\xi(s)ds+Ke^{-\frac{\d}4(T-t)}\dbE|\vartheta|^2.\ea\ee
Combining \eqref{wtYZ} and  \eqref{Eeta}, by $\eta_{\scT}(t)=\Sigma(t,\a(t))\wt\eta_{\scT}(t)$, we obtain that $\dbE|\eta_{\scT}(\cd)|^2$ is uniformly bounded on $[0,T]$, i.e.,
\bel{E|eta|}\dbE|\eta_{\scT}(t)|^2\les K\int_t^{T}
e^{-{\d\over4}(s-t)}\xi(s)ds+Ke^{-\frac{\d}4(T-t)}\dbE|\vartheta|^2.\ee
Next, from \eqref{wtYZ} and \rf{Eeta}, for $t\in[0,T]$, we have (see \rf{wteta} again)
\begin{align}
&\int_t^{T}e^{-{\d\over4}(s-t)}\dbE|\z_{\scT}(s)|^2ds\les
K\int_t^{T}e^{-{\d\over4}(s-t)}\dbE|\wt\z_{\scT}(s)|^2ds\nonumber\\
&\les Ke^{-\frac\d4(T-T_0-t)}\int_{T-T_0}^{T}e^{-{\d\over4}(s-(T-T_0))}\dbE|\wt \z_{\scT}(s)|^2ds+K\int_t^{T-T_0}e^{-{\d\over4}(s-t)}\dbE|\wt\z_{\scT}(s)|^2ds\nonumber\\
\ns\ds&\les Ke^{-\frac\d4(T-T_0-t)}\int_{T-T_0}^{T}e^{-{\d\over4}(s-(T-T_0))}\dbE\(|\wt\z_{\scT}(s)-\wt C^{\Th_T}(s,\a(s))\wt\eta_{\scT}(s)|^2+|\wt\eta_{\scT}(s)|^2\)ds\nonumber\\
&\q+K\int_t^{T-T_0}e^{-{\d\over4}(s-t)}\dbE\(|\wt\z_{\scT}(s)-\wt C^{\Th_T}(s,\a(s))\wt\eta_{\scT}(s)|^2+|\wt\eta_{\scT}(s)|^2\)ds\nonumber\\
&\les Ke^{-\frac\d4(T-T_0-t)}\(\int_{T-T_0}^{T}e^{-\frac{\d}4(s-(T-T_0))}\zeta(s)ds+\int_{T-T_0}^{T}\int_s^{T}e^{-{\d\over4}(r-s)}\xi(r)drds\)\nonumber\\
&\q+K\(\int_t^{T}
e^{-{\d\over4}(s-t)}\xi(s)ds+\int_t^{T-T_0}\int_s^{T}e^{-{\d\over4}(r-s)}\xi(r)drds\)+Ke^{-\frac{\d}4(T-t)}\dbE|\vartheta|^2\nonumber\\
&\les K\int_t^T
e^{-{\d\over4}(s-t)}\xi(s)ds+Ke^{-\frac{\d}4(T-t)}\dbE|\vartheta|^2.\label{E|zeta|}\end{align}
Likewise,
\bel{E|zeta^M|}\ba{ll}
\ns\ds\int_t^{T}e^{-{\d\over4}(s-t)}\dbE\big[\sum_{\imath\neq\jmath}\lambda_{\imath\jmath}|\z^M_{\scT}(s,\jmath)|^2{\bf 1}_{[\a(s)=\imath]}\big]ds\\
\ns\ds=\int_t^{T}e^{-{\d\over4}(s-t)}\dbE\Big[\sum_{\imath\neq \jmath}
|\wt\z^M_{\scT}(s,\jmath)+E(\imath,\jmath)\wt\eta_{\scT}(s)|^2
\l_{\imath\jmath}{\bf1}_{\{\a(s)=\imath\}}\]ds\\
\ns\ds\les K\int_t^T e^{-{\d\over4}(s-t)}\xi(s)ds+Ke^{-\frac{\d}4(T-t)}\dbE|\vartheta|^2.\ea\ee
Combining \rf{E|eta|}--\rf{E|zeta^M|}, we get \rf{stabsde} by taking $\vartheta=0$.

\it Step 5. Stability estimates. \rm For $T'>T\ges T_0$, we have (recall \rf{BSDET})
\bel{d[eta-eta]}\ba{ll}
\ns\ds d(\eta_{\scT}-\eta_{\scT'})=-\((A^{\Th_T}\eta_{\scT}+C^{\Th_T}\z_{\scT}
+\f_{\scT}-(A^{\Th_{T'}}\eta_{\scT'}+C^{\Th_{T'}}\z_{\scT'}+\f_{\scT'})
\)dt\\
\ns\ds\qq\qq\qq\qq+(\z_{\scT}-\z_{\scT'})dW+(\z_{\scT}^M-\z_{\scT'}^M)
dM\\
\ns\ds=-\((A^{\Th_{T'}})^\top(\eta_{\scT}-\eta_{\scT'})+C^{\Th_{T'}}
(\z_{\scT}-\z_{\scT'})\\
\ns\ds\qq\qq+(A^{\Th_T}-A^{\Th_{T'}})\eta_{\scT}+(C^{\Th_T}-C^{\Th_{T'}})
\z_{\scT}+\f_{\scT}-\f_{\scT'}\)ds\\
\ns\ds\qq\qq+(\z_{\scT}-\z_{\scT'})dW+(\z_{\scT}^M-\z_{\scT'}^M)dM\\
\ns\ds\equiv -\((A^{\Th_{T'}})^\top(\eta_{\scT}-\eta_{\scT'})+C^{\Th_{T'}}
(\z_{\scT}-\z_{\scT'})+\D_{\scT,T'}(s)\)dt\\
\ns\ds\qq\qq+(\z_{\scT}-\z_{\scT'})dW
+(\z_{\scT}^M-\z_{\scT'}^M)dM,\ea\ee
where
$$\ba{ll}
\ns\ds\D_{\scT,T'}(s)=(A^{\Th_T}-A^{\Th_{T'}})^\top\eta_{\scT}
+(C^{\Th_T}-C^{\Th_{T'}})^\top\z_{\scT}+\f_{\scT}-\f_{\scT'},\\
\ns\ds\f_{\scT }=P_{\scT}b+(C^{\Th_{\scT}})^\top P_{\scT}\si+q+\Th_{\scT}^\top r,\q\f_{\scT'}= P_{\scT'}b+(C^{\Th_{T'}})^\top P_{\scT'}\si+q+\Th_{\scT'}^\top r.\ea$$
Since $T'>T$, we have $P_{\scT}(t,\imath)\les P_{\scT'}(t,\imath)\les P_{\sc\i}(t,\imath)$, it follows that
$$0\les P_{\sc{T'}}(t,\imath)-P_{\scT}(t,\imath)\les P_{\sc\i}(t,\imath)-P_{\scT}(t,\imath)\les Ke^{-\d(T-t)}I.$$
For $t\in[0,T]$, we have
\bel{DeltaTT}\ba{ll}
\ns\ds\dbE|\D_{\scT,T'}(t)|^2\les
\dbE\(|(A^{\Th_T}-A^{\Th_{T'}})^\top\eta_{\scT}
+(C^{\Th_T}-C^{\Th_{T'}})^\top\z_{\scT}|
+|\f_{\scT}-\f_{\scT'}|\)^2\\
\ns\ds\les\dbE\(|B(\Th_{\scT}-\Th_{\scT'})\eta_{\scT}|
+|D(\Th_{\scT}-\Th_{\scT'})\z_{\scT}|\\
\ns\ds\qq+|(P_{\scT}-P_{\scT'})b|+|(C^{\Th_T}P_{\scT}
-C^{\Th_{T'}}P_{\scT'})\si|+|(\Th_{\scT}-\Th_{\scT'})^\top r|\)^2\\
\ns\ds\les Ke^{-2\d(T-t)}\dbE\(|b|+|\si|+|r|+|\eta_{\scT}|+|\z_{\scT}|\)^2\\
\ns\ds\les Ke^{-2\d(T-t)}\(\xi(t)+\dbE|\eta_{\scT}(t)|^2+\dbE|\z_{\scT}(t)|^2\).\ea\ee
Note that  \eqref{d[eta-eta]} is parallel with BSDE \eqref{BSDEwteta} with different non-homogeneous terms and terminal condition only. Similar to Steps 2--3,  it follows that
\begin{align}
&\dbE|\eta_{\scT}(t)-\eta_{\scT'}(t)|^2+\dbE\int_t^{T}
e^{-{\d\over4}(s-t)}|\z_{\scT}(s)-\z_{\scT'}(s)|^2ds\nonumber\\
&\qq+\dbE\int_t^{T}e^{-{\d\over4}(s-t)}\sum_{\jmath\neq\imath}\lambda_{\imath\jmath}|\z^M_{\scT}(s,\jmath)-\z^M_{\scT'}(s,\jmath)|^2{\bf 1}_{[\a(s)=\imath]}ds\nonumber\\
&\les K\int_t^{\scT} e^{-\frac\d4(s-t)}e^{-2\d(T-s)} \dbE|\D_{\scT,\scT'}|^2ds+Ke^{-{\d\over4}(T-t)}\dbE|\eta_{\scT}(T)-\eta_{\scT'}(T)|^2\nonumber\\
&\les K\int_t^{\scT} e^{-\frac\d4(s-t)}e^{-2\d(T-s)} \(\xi(s)+\dbE|\eta_{\scT}(s)|^2+\dbE|\z_{\scT}(s)|^2\)ds+Ke^{-{\d\over4}(T-t)}\dbE|\eta_{\scT}(T)-\eta_{\scT'}(T)|^2\nonumber\\
&\les e^{-\frac\d 8(T-t)}\int_t^{T'}e^{-\frac\d4(s-t)}\xi(s)ds.\label{deltaeta}
\end{align}
In the last step, we use $\eta_{\scT}(T)=0$ and \eqref{E|eta|} (taking $T=T'$ and $t=T$).

\end{proof}

\ms

\begin{proof}[Proof of Proposition \ref{festimateprop}]
(1) We consider the BSDE \eqref{BSDEwteta}. By \rf{deta} and \eqref{deta2}, we have for $t\in[0, T]$,
$$\ba{ll}
\ns\ds{d\over dt}\dbE|\wt\eta_{\scT}(t)|^2+{\d\over4}\dbE|\wt\eta_{\scT}(t)|^2\ges
\dbE|\wt\z_{\scT}(t)-\wt C^{\Th_T}(t,\a(t))\wt\eta_{\scT}(t)|^2\\
\ns\ds\q+\dbE\sum_{\imath\neq\jmath}
|\wt\z^M_{\scT}(t,\jmath)-E(\imath,\jmath)\wt\eta_{\scT}(t)|^2
\l_{\imath\jmath}{\bf1}_{\{\a(t)=\imath\}}-K\xi(t)-K\dbE|\wt\eta_{\scT}(t)|^2{\bf 1}_{t\in[T-T_0,T]}.\ea$$
Grownwall's inequality implies that
$$\ba{ll}
\ns\ds\dbE|\wt\eta_{\scT}(t)|^2-e^{-{\d\over4}t}|\wt\eta_{\scT}(0)|^2\ges\int_0^t
e^{-{\d\over4}(t-s)}\(\dbE|\wt\z_{\scT}(s)-\wt C^{\Th_T}(s,\a(s))\wt\eta_{\scT}(s)|^2\\
\ns\ds\q+\sum_{\imath\neq\jmath}
|\wt\z^M_{\scT}(s,\jmath)-E(\imath,\jmath)\wt\eta_{\scT}(s)|^2
\l_{\imath\jmath}{\bf1}_{\{\a(s)=\imath\}}-K\xi(s)+K\dbE|\wt\eta_{\scT}(s)|^2{\bf 1}_{s\in[T-T_0,T]}\)ds.\ea$$
Hence, for $t\in[0,T]$, (see \rf{wteta} again)
\bel{|z|}\ba{ll}
\ns\ds\int_0^te^{-{\d\over4}(t-s)}\dbE|\z_{\scT}(s)|^2ds\les
K\int_0^te^{-{\d\over4}(t-s)}\dbE|\wt\z_{\scT}(s)|^2ds\\
\ns\ds\les K\int_0^te^{-{\d\over4}(t-s)}\dbE\(|\wt\z_{\scT}(s)-\wt C^{\Th_T}(s,\a(s))\wt\eta_{\scT}(s)|^2+|\wt C^{\Th_T}(s,\a(s))\wt\eta_{\scT}(s)|^2\)dr\\
\ns\ds\les K\dbE\[|\wt\eta_{\scT}(t)|^2+\int_0^t
e^{-{\d\over4}(t-s)}\(|\wt\eta_{\scT}(s)|^2+\xi(s)\)ds+\int_0^te^{-{\d\over4}(t-s)}\dbE|\wt\eta_{\scT}(s)|^2{\bf 1}_{s\in[T-T_0,T]}ds\]\\
\ns\ds\les K\int_0^\i e^{-{\d\over4}|t-s|}\xi(s)ds+K e^{-\frac\d 4(T-t)}\dbE|\vartheta|^2.\ea\ee
%
%
Taking $\vartheta=0$,  we get \rf{festimate1}.

(2). Consider \eqref{d[eta-eta]} and  \eqref{DeltaTT}, and take $T'=\i$. By virtue of \eqref{|z|}, using \eqref{deltaeta}, we have
\bel{|z-2|}\ba{ll}
\ns\ds\int_0^te^{-{\d\over4}(t-s)}\dbE|\z_{\scT}(s)-\z_{\sc\i}(s)|^2ds\\
\ns\ds\les K\dbE\[|\wt\eta_{\scT}(t)-\wt\eta_{\sc\i}(t)|^2+\int_0^t
e^{-{\d\over4}(t-s)}\(|\wt\eta_{\scT}(s)-\wt\eta_{\sc\i}(s)|^2+|\D_{\sc T,\i}(s)|^2\)ds\\
\ns\ds\qq+\int_0^t\dbE|\wt\eta_{\scT}(s)-\wt\eta_{\sc\i}(s)|^2{\bf 1}_{s\in[T-T_0,T]}ds\]\\
\ns\ds\les Ke^{-\frac\d 8(T-t)}\int_0^\i e^{-{\d\over4}|t-s|}\xi(s)ds.\ea\ee

\ms

(3) Note that
$$\ba{ll}
\ns\ds\dbE|v_{\sc\i}(s,\a(s))-v_{\scT}(s,\a(s))|^2\\
\ns\ds=\dbE\big|\wt R_{\sc\i}(\a(s))^{-1}[D^\top P_{\sc\i}(\a(s))\si(s)+B^\top\eta_{\sc\i}(s)+D^\top
\z_{\sc\i}(s)+r(s)]\\
\ns\ds\qq-\wt R_{\scT}(\a(s))^{-1}[D^\top P_{\scT}(s,\a(s))\si(s)+B^\top\eta_{\scT}(s)+D^\top
\z_{\scT}(s)+r(s)]\big|^2\\
\ns\ds\les K\dbE\[\big|\wt R_{\sc\i}(\a(s))^{-1}\big|\(|P_{\sc\i}(\a(s))-P_{\scT}(s,\a(s))|\,
|\si(s)|+|\eta_{\sc\i}(s)-\eta_{\scT}(s)|+|\z_{\sc\i}(s)-\z_{\scT}(s)|\)\\
\ns\ds\q+\big|\wt R_{\sc\i}(\a(s))^{-1}-\wt R_{\scT}(s,\a(s))^{-1}\big|\Big(|\eta_{\sc\i}(s)|
+|\z_{\sc\i}(s)|+|\si(s)|+|r(s)|\)\]^2\\
\ns\ds\les K\dbE\[|P_{\sc\i}(\a(s))-P_{\scT}(s,\a(s))|^2|\si(s)|^2+|\eta_{\sc\i}(s)-\eta_{\scT}(s)|^2+|\z_{\sc\i}(s)-\z_{\scT}(s)|^2
\)\\
\ns\ds\q+\big|P_{\sc\i}(\a(s))-P_{\scT}(s,\a(s))\big|^2\(|\eta_{\sc\i}(s)|^2
+|\z_{\sc\i}(s)|^2+|\si(s)|^2+|r(s)|^2\)\]\ea$$
$$\ba{ll}
\ns\ds\les K\dbE\[\(|\eta_{\sc\i}(s)-\eta_{\scT}(s)|^2+|\z_{\sc\i}(s)
-\z_{\scT}(s)|^2~~~~~~~~~~~~~~~~~~~~~~~~~~~~~~~~~~~~~~~~~~~~~~~~~~~~~~~~\\
\ns\ds\q+e^{-2\d(T-s)}\(|\eta_{\sc\i}(s)|^2
+|\z_{\sc\i}(s)|^2+|\si(s)|^2+|r(s)|^2\)\].\ea$$
Hence, for $t\in[0,T]$, we have
\bel{v-v1}\ba{ll}
\ns\ds\int_0^te^{-\frac\d4(t-s)}\dbE|v_{\sc\i}(s,\a(s))-v_{\scT}(s,\a(s))|^2 ds\\
\ns\ds\les K\dbE\int_0^te^{-\frac\d4(t-s)}\(|\eta_{\sc\i}(s)-\eta_{\scT}(s)|^2
+|\z_{\sc\i}(s)-\z_{\scT}(s)|^2\\
\ns\ds\qq\qq\qq+e^{-2\d(T-s)}(|\eta_{\sc\i}(s)|^2+|\z_{\sc\i}(s)|^2+|\si(s)|^2
+|r(s)|^2)\)ds\\
\ns\ds\les Ke^{-{\d\over8}(T-t)}\(\int_0^\i e^{-{\d\over4}|t-s|}\xi(s)ds\).\ea\ee
(4) By \eqref{deta} and \eqref{deta2}, we have
$$\ba{ll}
\ns\ds{d\over dt}\dbE|\wt\eta_{\scT}(t)|^2\ges
\dbE|\wt\z_{\scT}(t)-\wt C^{\Th_T}(t,\a(t))\wt\eta_{\scT}(t)|^2-K\xi(t)-K\dbE|\wt\eta_{\scT}(t)|^2{\bf 1}_{t\in[T-T_0,T]}.\ea$$
Integrating both sides on $[0,T]$ implies that
$$\ba{ll}
\ns\ds\dbE|\wt\eta_{\scT}(T)|^2-|\wt\eta_{\scT}(0)|^2\\
\ns\ds\ges\int_0^T
\dbE|\wt\z_{\scT}(s)-\wt C^{\Th_T}(s,\a(s))\wt\eta_{\scT}(s)|^2-K\xi(s)-K\dbE|\wt\eta_{\scT}(s)|^2{\bf 1}_{s\in[T-T_0,T]}\)ds.\ea$$
Hence, taking $\vartheta=0$,
\begin{align*}\int_0^T
&\dbE|\z_{\scT}(s)|^2ds\les K\int_0^T\dbE|\wt\z_{\scT}(s)|^2ds\\
&\les
K\int_0^T
\dbE|\wt\z_{\scT}(s)-\wt C^{\Th_T}(s,\a(s))\wt\eta_{\scT}(s)|^2ds+K\int_0^T
\dbE|\wt\eta_{\scT}(s)|^2ds\\
&\les K\int_0^T
\xi(s)ds+K\int_0^T\dbE|\wt\eta_{\scT}(s)|^2ds+\dbE|\wt\eta_{\scT}(T)|^2\\
&\les K\int_0^T
\xi(s)ds+K (T+1)\sup_{s\ges 0}\int_{0}^\infty e^{-\frac\d 4|s-r|}\xi(r)dr.
\end{align*}

Similarly, one has
$$\ba{ll}
\ns\ds{d\over dt}\dbE|\wt\eta_{\sc\i}(t)|^2\ges
\dbE|\wt\z_{\sc\i}(t)-\wt C^{\Th_{\sc\i}}(t,\a(t))\wt\eta_{\sc\i}(t)|^2-K\xi(t).\ea$$
Therefore, we have
\begin{align*}\int_0^T
&\dbE|\z_{\sc\i}(s)|^2ds\les K\int_0^T\dbE|\wt\z_{\sc\i}(s)|^2ds\\
&\les
K\int_0^T
\dbE|\wt\z_{\sc\i}(s)-\wt C^{\Th_{\sc\i}}(s,\a(s))\wt\eta_{\scT}(s)|^2ds+K\int_0^T
\dbE|\wt\eta_{\sc\i}(s)|^2ds\\
&\les K\int_0^T
\xi(s)ds+K\int_0^T\dbE|\wt\eta_{\sc\i}(s)|^2ds+\dbE|\wt\eta_{\sc\i}(T)|^2\\
&\les K\int_0^T
\xi(s)ds+K (T+1)\sup_{s\ges 0}\int_{0}^\infty e^{-\frac\d 4|s-r|}\xi(r)dr.
\end{align*}

(5). By \eqref{d[eta-eta]} and \eqref{DeltaTT} (letting $T'=\i$), we have
$$\ba{ll}
\ns\ds{d\over dt}\dbE|\wt\eta_{\scT}(t)-\wt\eta_{\sc\i}(t)|^2\ges
\dbE|\wt\z_{\scT}(t)-\wt C^{\Th_{\scT}}(t,\a(t))\wt\eta_{\scT}(t)-\wt\z_{\sc\i}(t)+\wt C^{\Th_{\sc\i}}(t,\a(t))\wt\eta_{\sc\i}(t)|^2-K\dbE|\D_{\scT,\sc\i}(t)|^2.\ea$$
Thus, it follows that
\begin{align*}\int_0^T
&\dbE|\z_{\scT}(s)-\z_{\sc\i}(s)|^2ds\les K\int_0^T\dbE|\wt\z_{\scT}(s)-\wt\z_{\sc\i}(s)|^2ds\\
&\les
K\int_0^T
\dbE|\wt\z_{\scT}(s)-\wt C^{\Th_{\scT}}(s,\a(s))\wt\eta_{\scT}(t)-\wt\z_{\scT}(s)+\wt C^{\Th_{\sc\i}}(s,\a(s))\wt\eta_{\sc\i}(s)|^2ds\\
&\q+K\int_0^T
\dbE|\wt\eta_{\scT}(s)-\wt\eta_{\sc\i}(s)|^2ds+\dbE|\wt\eta_{\scT}(T)-\wt\eta_{\sc\i}(T)|^2\\
&\les K\int_0^T
\dbE|\D_{\scT,\sc\i}(t)|^2ds+K\int_0^Te^{-\frac\d 2(T-s)}\int_0^\i  e^{-{\d\over4}|r-s|}\xi(r)drds+\dbE|\wt\eta_{\sc\i}(T)|^2\\
&\les K\int_0^T
e^{-2\d(T-t)}\(\xi(t)+\dbE|\eta_{\scT}(t)|^2+\dbE|\z_{\scT}(t)|^2\)ds+K \sup_{s\ges 0}\int_{0}^\infty e^{-\frac\d 4|s-r|}\xi(r)dr\\
&\les K\sup_{s\ges 0}\int_{0}^\infty e^{-\frac\d 4|s-r|}\xi(r)dr
\end{align*}
In the last step, we used \eqref{festimate1} (with $t=T$) for $\zeta_{\scT}(\cd)$.

(6) We suppress the index $(x,\imath)$ in this part of proof.  Note that
\bel{closed1}\left\{\2n\ba{ll}
\ds d\bar X_{\scT}(t)=[A^{\Th_T}(t,\a(t))\bar X_{\scT}(t)+B(\a(t))v_{\scT}(t,\a(t))+b(t)]dt\\
\ns\ds\qq\qq\q+[C^{\Th_T}(t,\a(t))\bar X_{\scT}(t)+D(\a(t)) v_{\scT}(t,\a(t))+\si(t)]dW(t),\\
\ns\ds\bar X_{\scT}(0)=x,\ea\right.\ee
with $A^{\Th_T}(\cd\,,\a(\cd))$ and $C^{\Th_T}(\cd\,,\a(\cd))$ being given by \rf{A^Th_T}. Now, let $\Si_{\sc\i}(\cd)\in\hSi$ satisfy \rf{Ly1}. Applying It\^o's formula to the map $t\mapsto\lan\Si_{\sc\i}(\a(t))\bar X_{\scT}(t),\bar X_{\scT}(t)\ran$, we have, (see \rf{Ly2}) suppressing $\a(t)$,
$$\ba{ll}
\ns\ds{d\over dt}\dbE\lan\Si_{\sc\i}(\a(t))\bar X_{\scT}(t),\bar X_{\scT}(t)\ran\\
\ns\ds=\dbE\Big\lan\(\L[\Si_{\sc\i}]+\Si_{\sc\i}A^{\Th_{\scT}}
+(A^{\Th_{\scT}})^\top\Si_{\sc\i}+(C^{\Th_{\scT}})^\top\Si_{\sc\i}
C^{\Th_{\scT}}\)\bar X_{\scT}(t),\bar X_{\scT}(t)\Big\ran\\
\ns\ds\q+2\dbE\lan\Si_{\sc\i}(Bv_{\scT}(t)+b(t)),\bar X(t)\ran+2\dbE\lan\Si_{\sc\i}C^{\Th_{\scT}}\bar X_{\scT}(t),Dv_{\scT}(t,\a(t))+\si(t)\ran\\
\ns\ds\q+\dbE\lan\Si_{\sc\i}(Dv_{\scT}(t,\a(t))+\si(t),Dv_{\scT}(t,\a(t))+\si(t)\ran\\
\ns\ds\les-{\d\over2}\dbE\lan\Si_{\sc\i}(\a(t))\bar X_{\scT}(t),\bar X_{\scT}(t)\ran+K\dbE\(|b(t)|^2+|\si(t)|^2+|v_{\scT}(t,\a(t))|^2\).\ea$$
Note that $$\dbE|v_{\scT}(t,\a(t))|^2\leq K\(\xi(t)+\dbE|\eta_{\scT}(t)|^2+\dbE|\zeta_{\scT}(t)|^2\).$$
By Grownwall's inequality, we have
$$\ba{ll}
\ns\ds\dbE\lan\Si_{\sc\i}(\a(t))\bar X_{\scT}(t),\bar X_{\scT}(t)\ran\\
\ns\ds\les K e^{-{\d\over2}t}|x|^2+K\int_0^t e^{-{\d\over2}(t-s)}\dbE\Big(|b(s)|^2+|\si(s)|^2+|v_{\scT}(s,\a(s))|^2\Big)ds\\
\ns\ds\les K e^{-{\d\over2}t}|x|^2+K\int_0^\i e^{-{\d\over4}|t-s|}\xi(s)ds,\qq 0\les t<T.\ea$$
Therefore, \eqref{boundEXTX} holds for $\bar X_{\scT}(\cd)$. The proof for $\bar X_{\sc\i}(\cd)$ is identical.
\end{proof}


\begin{thebibliography}{90}
\addtolength{\itemsep}{-1.0ex}

\bibitem{Basei-Pham-2019} M.~Basei and H.~Pham, \it A weak martingale approach to linear-quadratic McKean-Vlasov stochastic control problems, \sl J. Optim. Theory Appl., \rm 181 (2019), 347--382.

\bibitem{Bayraktar-Jian-2025} E.~Bayraktar and J.~Jian, \it Ergodicity and turnpike properties of linear-quadratic mean field control problems, \rm arXiv preprint arXiv:2502.08935.

\bibitem{Bian-Zheng-2019} B.~Bian and H.~Zheng, \it Turnpike property and convergence rate for an investment and consumption model, \sl Math Financial Economics, \rm 13 (2019), 227--251.


\bibitem{Breiten-Pfeiffer-2020} T.~Breiten and L.~Pfeiffer, \it On the turnpike property and the receding-horizon method for linear-quadratic optimal control problems, \sl SIAM J. Control Optim., \rm 58 (2020), 1077--1102.

\bibitem{Carlson-Haurie-Leizarowitz-1991} D.~A.~Carlson, A.~B.~Haurie, and A.~Leizarowitz, \sl Infinite Horizon Optimal Control --- Deterministic and Stochastic Systems, 2nd ed., \rm Springer-Verlag, Berlin, 1991.

\bibitem{Chen-Luo-2023} Y.~Chen and P.~Luo, \it Turnpike properties for stochastic backward linear-quadratic optimal problems, \rm
    arXiv preprint arXiv:2309.03456.

\bibitem{Conforti-2023} G.~Conforti, \it Coupling by reflection for controlled diffusion processes: Turnpike property and large time behavior of Hamilton-Jacobi-Bellman equations, \sl Ann. Appl. Probab., \rm 33 (2023), 4608--4644.

\bibitem{Cox-Huang-1992} J.~C.~Cox and C.-F.~Huang, \it A continuous-time portfolio turnpike theorem, \sl J. Economic Dynamics \& Control, \rm 16 (1992), 491--507.

\bibitem{Damm-Grune-Stieler-Worthmann-2014} T.~Damm, L.~Gr\"{u}ne, M.~Stieler, and K.~Worthmann, \it An exponential turnpike theorem for dissipative discrete time optimal control problems, \sl SIAM J. Control Optim., \rm 52 (2014), 1935--1957.

\bibitem{Dorfman-Samuelson-Solow-1958} R.~Dorfman, P.~A.~Samuelson, and R.~M.~Solow, \sl Linear Programming and Economics Analysis, \rm McGraw-Hill, New York, 1958.

\bibitem{Dybvig-Rogers-Back-1999} P.~H.~Dybvig, L.~C.~G.~Rogers, and K.~Back, \it Porfolio turnpikes, \sl The Review Financial Studies, \rm 12 (1999), 165--195.

\bibitem{Faulwasser-Grune-2022} T.~Faulwasser and L.~Gr\"{u}ne, \it Turnpike properties in optimal control: An overview of discrete-time and continuous-time results, \sl Numerical Control: Part A, \rm 23 (2022), 367--400.

\bibitem{Geng-Zariphopoulou-2025} T.~Geng and T.~Zariphopoulou, \it Temporal and sptial turnpikes in Iot-diffusion markets under forward performance criteria, \sl Numer. Alg. Control. Optim., \rm 15 (2025), 243--272.

\bibitem{Grune-Guglielmi-2018} L.~Gr\"{u}ne and R.~Guglielmi, \it Turnpike properties and strict dissipativity for discrete time linear quadratic optimal control problems, \sl SIAM J. Control Optim., \rm 56 (2018), 1282--1302.

\bibitem{Grune-Guglielmi-2021} L.~Gr\"{u}ne and R.~Guglielmi, \it On the relation between turnpike properties and dissipativity for continuous time linear quadratic optimal control problems, \sl Math. Control Relat. Fields, \rm 11 (2021), 169--188.

\bibitem{Huang-Zariphopoulou-1999} C.-F.~Huang and T.~Zariphopoulou, \it Turnpike behavior of long-term invstments, \sl Finance Stochast., \rm 3 (1999), 15--34.

\bibitem{Huang-Li-Yong-2015} J.~Huang, X.~Li, and J.~Yong, \it A linear-quadratic optimal control problem for mean-field stochastic differential equations in infinite horizon, \sl Math. Control Relat. Fields, \rm 5 (2015), 97--139.

\bibitem{Humberman-Ross-1983} G.~Humberman and S.~Ross, \it Portfolo turnpike theorems, risk aversion, and regularly varying utlity functions, \sl Econometrica, \rm 51 (1983), 1345--1361.



\bibitem{Jin-1998} X.~Jin, \it Consumption and porfolio turnpike theorems in a continuous-time finance model, \sl J. Economic Dynamics \& Control, \rm 22 (1998), 1001--1026.


\bibitem{Jian-Jin-Song-Yong-2024} J.~Jian, S.~Jin, Q.~Song, and J.~Yong, \it Long-Time Behaviors of Stochastic Linear-quadratic Optimal Control Problems, \rm arXiv preprint arXiv:2409.11633.

\bibitem{Leland-1972} H.~Leland, \it On turnpike portfolios, \sl Mathematical Models in Investment and finance, G.~P.~Szego, Karl Shell eds, \rm 6--15.


\bibitem{Lu-2020} Q.~L\"{u}, \it Stochastic linear quadratic optimal control problems for mean-field stochastic evolution equation, \sl ESAIM Control Optim. Calc. Var., \rm 26 (2020), 127.

\bibitem{Lou-Wang-2019} H.~Lou and W.~Wang, \it Turnpike properties of optimal relaxed control problems, \sl ESAIM Control Optim. Calc. Var., \rm 25 (2019), 74.


\bibitem{Mao-1999} X.~Mao, \it Stability of stochastic differential equations with Markovian switching, \sl Stoch. Proc. Appl.,\rm 79 (1999) 45--67.

\bibitem{Marimon-1989} R.~Marimon, \it Stochastic turnpike property and stationary equilibrium, \sl J. Economic Theory, \rm 47 (1989), 282--306.

\bibitem{McKenzie-1976} L.~W.~McKenzie, \it Turnpike theory, \sl Econometrica, \rm 44 (1976), 841--865.

\bibitem{Mei-Wang-Yong-2025} H.~Mei, R. Wang, and J.~Yong, \it Turnpike property of stochastic linear-quadratic optimal control problems in large horizons with regime switching I: Homogeneous cases.\rm preprint.

\bibitem{Mei-Wei-Yong-2024}H. Mei, Q. Wei, and J. Yong, \it Linear-quadratic optimal control problem for mean-field stochastic differential equations with a type of random coefficients, \rm Numer. Algebra, Control Optim., 14 (2024), 813--852.

\bibitem{Mei-Wei-Yong-2025} H.~Mei, Q.~Wei, and J.~Yong, \it Linear-quadrartic optimal control for mean-field  stochastic differential equations in infinite-horizon with regime switching, \sl Chin. Ann. Math. \rm to appear.

\bibitem{Neumann-1945} J.~von Neumann, \it A model of general economic equilibrium, \sl Rev. Econ. Stud., \rm 13 (1945), 1--9.

\bibitem{Porretta-Zuazua-2013} A.~Porretta and E.~Zuazua, \it Long time versus steady state optimal control, \sl SIAM J. Control Optim., \rm 51 (2013), 4242--4273.


\bibitem{Ramsey-1928} F.~P.~Ramsey, \it A mathematical theory of saving, \sl The economic journal, \rm 38 (1928), 543--559.

\bibitem{Ross-1974} S.~A.~Ross, \it Portfolio turnpike theorems for constant policies, \sl J. Financial Economics, \rm 1 (1974), 171--198.

\bibitem{Sakamoto-Zuazua-2021} N.~Sakamoto and E.~Zuazua, \it The turnpike property in nonlinear optimal control---A geometric approach, \sl Automatica, \rm 134 (2021), 109939.

\bibitem{Schiessl-Baumann-Faulwasser-Grune-2024} \rm J.~Schiessl, M.~H.,~Baumann, T.~Faulwasser, and L.~Gr\"une, L., \it On the relationship between stochastic turnpike and dissipativity notions,  \rm arXiv:2311.07281v2 [math.OC] 21 Aug 2024.

\bibitem{Skorokhod-1989} A.~V.~Skorokhod, \sl Asymptotic Methods in the Theory of Stochastic Differential Equations, \rm AMS, Providence, RI., 1989.

\bibitem{Sun-2017} \rm J.~Sun, \it Mean-field stochastic linear quadratic optimal control problems: Open-loop solvabilities, \sl ESAIM Control Optim. Calc. Var., \rm 23 (2017), 1099--1127.


\bibitem{Sun-Wang-Yong-2022} \rm J.~Sun, H.~Wang, and J.~Yong, \it Turnpike properties for stochastic linear-quadratic optimal control problems, \sl Chin. Ann. Math., \rm Ser B, 43 (2022), 999--1022.

\bibitem{Sun-Yong-2020} \rm J.~Sun and J.~Yong, \sl Stochastic Linear-quadratic Optimal Control Theory: Open-Loop and Closed-Loop Solutions, \rm Springer Briefs Math., Springer, 2020.

\bibitem{Sun-Yong-2024S}  \rm J.~Sun and J.~Yong, \it Turnpike properties for mean-field linear-quadratic optimal control problems, \sl SIAM J. Control Optim., \rm 62 (2024), 752--775.

\bibitem{Sun-Yong-2024J} \rm J.~Sun and J.~Yong, \it Turnpike properties for stochastic linear-quadratic optimal control problems with periodic coefficients, \sl J. Diff. Equ., \rm 400 (2024), 189--229.

\bibitem{Trelat-Zuazua-2015} E.~Tr\'{e}lat and E.~Zuazua, \it The turnpike property in finite-dimensional nonlinear optimal control, \sl J. Diff. Equ., \rm 258 (2015), 81--114.

\bibitem{Yin-Zhang-2010} G.~G.~Yin and Q.~Zhang, \sl Continuous-Time Markov Chains and Applications --- A Two-Time-scale Approach, Second eds., \rm Springer, 2010.

\bibitem{Yin-Zhu-2010} G.~G.~Yin and C.~Zhu, \sl Hybrid Switching Diffusions --- Properies and Applications, \rm Springer, 2010.
\bibitem{Yong-2013} J.~Yong, \it Linear-quadratic optimal control problems for mean-field stochastic differential equations, \sl SIAM J. Control Optim., \rm 51 (2013), 2809--2838.


\bibitem{Yong-2018} J.~Yong, \sl Optimization Theory: A Concise Introduction, \rm World Scientific, Singapore, 2018.


\bibitem{Zaslavski-2006} A.~J.~Zaslavski, \sl Turnpike Properties in the Calculus of Variations and Optimal Control, \rm Nonconvex Optim. Appl. 80, Springer, New York, 2006.




\bibitem{Zaslavski-2019} A.~J.~Zaslavski, \sl Turnpike Theory of Continuous-Time Linear Optimal Control Problems, \rm Springer Optim. Appl. 148, Springer, 2019.


\bibitem{Zuazua-2017} E.~Zuazua, \it Large time control and turnpike properties for wave equations, \rm Ann. Rev. Control, \rm 44 (2017), 199--210.

\bibitem{Zhang-Li-Xiong-2021} X.~Zhang, X.~Li, and J.~Xiong, \it Open-loop and closed-loop solvabilities for stochastic linear quadratic optimal control problems of Markovian regime switching system, \sl ESAIM: Control, Optim. Calc. Var., \rm 27(2021), 69.

\end{thebibliography}
\end{document}